\documentstyle[12pt]{article}
\newcommand{\sect}[1]{\section{#1}\setcounter{equation}{0}}

\font\mbn=msbm10 scaled \magstep1
\font\mbs=msbm7 scaled \magstep1
\font\mbss=msbm5 scaled \magstep1
\newfam\mbff
\textfont\mbff=\mbn
\scriptfont\mbff=\mbs
\scriptscriptfont\mbff=\mbss\def\mbf{\fam\mbff}
\def\Re{{\mbf R}}

\def\Z{{\mbf Z}}
\def\Co{{\mbf C}}
\def\To{{\mbf T}}

\def\Di{{\mbf D}}
\newtheorem{Th}{Theorem}[section]
\newtheorem{Lm}[Th]{Lemma}
\newtheorem{C}[Th]{Corollary}
\newtheorem{D}[Th]{Definition}
\newtheorem{Proposition}[Th]{Proposition}
\newtheorem{R}[Th]{Remark}
\author{Alexander Brudnyi\thanks{Research supported in part by NSERC.
}\\
Department of Mathematics\\
Ben Gurion University of the Negev, Beer-Sheva\\
Israel}
\title{TOPOLOGY OF THE MAXIMAL IDEAL SPACE OF $H^{\infty}$}
\date{May 5, 1999} 
\begin{document} 
%==================================
%==================================
\maketitle
\begin{abstract}
{We study the structure of the maximal ideal space 
$M(H^{\infty})$ of the 
algebra $H^{\infty}=H^{\infty}(\Di)$ of bounded analytic functions defined
on the open unit disk $\Di\subset\Co$. 
Based on the fact that $dim\ M(H^{\infty})=2$ we prove for
$H^{\infty}$ the matrix-valued corona theorem.
Our results heavily rely on the 
topological construction describing maximal ideal spaces of certain
algebras of continuous functions defined on the covering spaces of compact 
manifolds.
}
\end{abstract}
\sect{\hspace*{-1em}. Introduction.}
Let $H^{\infty}$ be the uniform algebra of bounded analytic 
functions on the unit disk $\Di\subset\Co$ with the pointwise multiplication 
and with the norm
$$
||f||=\sup_{z\in\Di}|f(z)|.
$$
The maximal ideal space of $H^{\infty}$ is defined by
$$
M(H^{\infty})=\{\phi\ : \phi\in Hom(H^{\infty},\Co),\ \phi\neq 0\}
$$
equipped with the weak $*$ topology induced by the dual space of $H^{\infty}$.
It is a compact Hausdorff space. A function $f\in H^{\infty}$ can be 
thought of as a continuous function on $M(H^{\infty})$ via the Gelfand
transform $\hat f(\phi)=\phi(f)$\ $(\phi\in M(H^{\infty}))$. Evaluation at a
point of $\Di$ is an element of $M(H^{\infty})$, so $\Di$ is naturally
embedded into $M(H^{\infty})$, and $\hat f$ is an extension of $f$ 
to $M(H^{\infty})$. In what follows we avoid writing the hat for the Gelfand
transform of $f$.

In 1962, Carleson [C] proved the following famous corona theorem:

Let $f_{1},f_{2},...,f_{n}$ be functions in $H^{\infty}$ such that
$$
|f_{1}(z)|+|f_{2}(z)|+...+|f_{n}(z)|\geq\delta>0\ \ \ \ \ {\rm for\ all}\ 
z\in\Di.
$$
Then there are $H^{\infty}$ functions $g_{1},g_{2},...,g_{n}$ so that
$$
f_{1}(z)g_{1}(z)+f_{2}(z)g_{2}(z)+...+f_{n}(z)g_{n}(z)=1.
$$
This result is equivalent to the statement that $\Di$ is dense in 
$M(H^{\infty})$. In this paper we consider the following
matrix-valued version of the corona theorem:

Suppose that $A$ is a complex, commutative Banach algebra with unit and that
$\nu:H^{\infty}\longrightarrow A$ is a Banach algebraic morphism with dense
image. Let $a=(a_{ij})$ be a $(k\times n)$-matrix, $k\leq n$, with entries
in $A$.
\begin{Th}\label{new}
Assume that the minors of $a$ of order $k$ do not belong together to any
maximal ideal of $A$.  Then there exists a unimodular $(n\times n)$-matrix
$\widetilde a=(\widetilde a_{ij})$,  $\widetilde a_{ij}\in A$ 
for all $i,j$, so that $\widetilde a_{ij}=a_{ij}$ 
for $1\leq i\leq k$. 
\end{Th}
In the proof we use the remarkable theorem of Su\'{a}rez 
([S]) asserting that topological dimension of $M(H^{\infty})$ equals 2.
Similar {\em complement problems} have been posed and solved by V.Lin for
matrices with entries from various functional algebras (see [L] Th.3 and 4).\\
{\bf Examples.} (1) Let $A=H^{\infty}$ and $\nu=id$. According to the
corona theorem conditions of Theorem \ref{new} can be written as
\begin{equation}\label{matr}
\sum_{s\in S}|A_{s}(z)|\geq\delta>0\ \ \ \ \ {\rm for\ all}\ z\in\Di,
\end{equation}
where $\{A_{s}\}_{s\in S}$ is the family of minors of $a$ of order $k$. Theorem
\ref{new} then implies that there is an invertible $(n\times n)$-matrix with
entries in $H^{\infty}$ which completes $a$. This statement can be obtained
also as the combination of Fuhrmann-Vasyunin theorem and Tolokonnikov's lemma
(see [Ni] p.293 for both these results). Unlike the proof in [Ni] which heavily
relies upon the structure of $H^{\infty}$ in the disk, our proof of Theorem
\ref{new} is of a local topological nature. This allows to obtain similar
matrix corona theorems for a wide class of planar domains (see [Bru]).\\
(2) Let $f_{1},...,f_{d}$ be functions in $H^{\infty}$ such that the set
$$
S=\{\xi\in M(H^{\infty})\ :\ \max_{1\leq i\leq d}|f_{i}(\xi)|\leq 1\}
$$
is not empty. Let $A$ be closure in $C(S)$ of the restriction of
$H^{\infty}$ to $S$ and $\nu:H^{\infty}\longrightarrow A$ be the restriction
homomorphism. Then Theorem \ref{new} is applicable to matrix-functions
defined on $S$ with entries in $A$ whose minors satisfy (\ref{matr}) at each
point of $S$. (It is easy to prove that the maximal ideal space of $A$ is $S$.)

The second part of our paper deals with
the topological description of the analytical
part of $M(H^{\infty})$. Our results generalize well-known 
theorems of Hoffman ([Ho]).  
Based on our main construction we solve the complement problem 
for some uniform subalgebras of $H^{\infty}$ that should be thought of
as a generalization of the disk algebra of analytic functions on $\Di$ 
continuous up to the boundary. We proceed to the formulation of other 
results.

A pseudohyperbolic metric on $\Di$ is defined by
$$
\rho(z,w)=\left|\frac{z-w}{1-\overline{w}z}\right|\ \ \ \ \
(z,w\in\Di).
$$
For $x,y\in M(H^{\infty})$, the formula
$$
\rho(x,y)=\sup\{|f(y)| :\ f\in H^{\infty}, f(x)=0, ||f||\leq 1\}
$$
gives an extension of $\rho$ to $M(H^{\infty})$. The Gleason part of
$x\in M(H^{\infty})$ is then defined as $P(x)=\{y\in M(H^{\infty}) :\
\rho(x,y)<1\}$. It is well known that for $x,y\in M(H^{\infty})$ we
have $P(x)=P(y)$ or $P(x)\cap P(y)=\emptyset$. Hoffman's classification of
Gleason parts ([Ho]) shows that there are only two cases: either
$P(x)=\{x\}$\ $(x\in M(H^{\infty}))$ or $P(x)$ is an analytic disk.
The former case means that there is a continuous one-to-one and onto map
$L_{x}:\Di\longrightarrow P(x)$ such that $f\circ L_{x}\in H^{\infty}$
for every $f\in H^{\infty}$. Moreover, any analytic disk is contained
in a Gleason part, and any maximal analytic
disk is a Gleason part. Let $M_{a}\subset M(H^{\infty})$ denote the set
containing all non-trivial Gleason parts ($M_{a}$ will be called the
analytical part of $M(H^{\infty})$). In his work [S1] Su\'{a}rez proved that 
the set $M(H^{\infty})\setminus M_{a}$ of one point Gleason parts is totally
disconnected. One of the main goals of this paper is to describe the 
topology of $M_{a}$.

Let $l^{\infty}(G)$ be the algebra of bounded complex-valued functions on
a discrete group $G$ with the pointwise multiplication and the norm
$||x||=\sup_{g\in G}|x(g)|$. Let $\beta G$ be the {\em Stone-\v{C}ech
compactification} of $G$, i.e., the maximal ideal space of $l^{\infty}(G)$
equipped with the Gelfand topology. Then $G$ is naturally embedded into
$\beta G$ as an open everywhere dense subset, and topology on $G$ induced
by this embedding coincides with the original one, i.e., is discrete. Every 
function
$f\in l^{\infty}(G)$ has a unique extension $\hat f\in C(\beta G)$. Further,
the natural right and left actions of $G$ on itself are extendible to
continuous right and left actions of $G$ on $\beta G$.\\
Assume now that $G$ is the fundamental group of a compact Riemann surface
$M$ of genus $g\geq 2$.
\begin{Th}\label{class}
(a)\ There is an open everywhere dense subset $G_{in}\subset\beta G$
invariant with respect to the right and left actions of $G$ and a locally
trivial fiber bundle $p:E(M,G_{in})\longrightarrow M$ over $M$ with the
fiber $G_{in}$ such that $E(M,G_{in})$ is homeomorphic to $M_{a}$.\\
(b)\ Every bounded uniformly continuous (respectively, Lipschitz) with respect 
to the metric 
$\rho$ function $f$ admits a continuous extension $\hat f$ to
$M_{a}$  so that for every $x\in M_{a}$ the function $\hat f\circ L_{x}$ is 
uniformly continuous (respectively, Lipschitz) with respect to $\rho$ on 
$\Di$.\\
Here $L_{x}:\Di\longrightarrow P(x)$ determines the Gleason part $P(x)$
(see definition above).\\
(c)\ Let $F:Y\longrightarrow M_{a}$ be a continuous mapping of a connected
and locally connected space $Y$ to $M_{a}$. Then there is $x\in M_{a}$ so that
$F(Y)\subset P(x)$.
\end{Th}
In the next section we state our results needed in the proof of
Theorem \ref{class} and, also, some of consequences of this theorem.
%==========================================
\sect{\hspace*{-1em}. Maximal Ideal Spaces of Certain Banach Algebras.}
{\bf 2.1.} The proof of Theorem \ref{class} is based on the description of
maximal ideal spaces of some algebras of continuous functions
defined on regular coverings over compact manifolds.

Let $p:M_{G}\longrightarrow M$ be a Galois covering over a compact
manifold $M$ with a covering group $G$ and $d$ be a metric
on $M_{G}$ compatible with its topology and invariant with respect to the 
action of the group $G$. Denote by $C_{d}(M_{G})$ the algebra of bounded
functions on $M_{G}$ uniformly continuous with respect to $d$ equipped 
with the supremum norm. We describe the maximal ideal space of
$C_{d}(M_{G})$. 

In what follows we assume that $M$ has one of the following structures: 
continuous, Lipschitz, $C^{\infty}$, real analytic, complex. Then 
$M_{G}$ can be equipped with the same structure as $M$. Moreover, 
the group $G$ acts discretely on $M_{G}$ by mappings which preserve the
structure. 

It is well-known that $M_{G}$ can be thought of as a 
principle fibre bundle over $M$ with the discrete fibre $G$. Namely, there
is an open (finite) covering ${\cal U}=\{U_{i}\}_{i\in I}$ of $M$ by sets
homeomorphic to open Euclidean balls in some $\Re^{n}$ and a locally constant
cocycle $g=\{g_{ij}\}\in Z^{1}({\cal U}, G)$ 
such that $M_{G}$ is isomorphic (in the corresponding category) 
to the quotient space of the disjoint union 
$V=\sqcup_{i\in I}U_{i}\times G$ by the equivalence relation:
$U_{i}\times G\ni (x,g)\sim (x,gg_{ij})\in U_{j}\times G$.
The identification space is a fibre bundle with projection 
$p:M_{G}\longrightarrow M$ induced by the 
product projection $U_{i}\times G\longrightarrow U_{i}$, see, e.g. [Hi], Ch.1.
Let $\beta G$ be the {\em Stone-\v{C}ech compactification} of $G$ and
$E(M,\beta G)$ be a locally trivial fibre bundle over $M$ with the fibre
$\beta G$ constructed by the cocycle $g$. Then as above
$E(M,\beta G)$ is homeomorphic to the quotient space of the disjoint union 
$\widetilde V=\sqcup_{i\in I}U_{i}\times \beta G$ by
the equivalence relation:
$U_{i}\times\beta G\ni (x,\xi)\sim (x,\xi g_{ij})\in U_{j}\times\beta G$.
The projection $\widetilde p: E(M,\beta G)\longrightarrow M$ is induced by
the product projection $U_{i}\times\beta G\longrightarrow U_{i}$.
Note that there is a natural embedding $V\hookrightarrow \widetilde V$ 
induced by the embedding $G\hookrightarrow\beta G$. This embedding
commutes with the corresponding equivalence relations and so determines
an embedding of $M_{G}$ into $E(M,\beta G)$ as an open everywhere dense
subset. 
\begin{Th}\label{te1}
For every $f\in C_{d}(M_{G})$ there is a unique continuous extension
$\widehat f\in C(E(M,\beta G))$. The extended algebra coincides with
$C(E(M,\beta G))$. In particular, 
$E(M,\beta G)$ is homeomorphic to the maximal ideal space of 
$C_{d}(M_{G})$.
\end{Th}
We now describe some topological properties of $E(M,\beta G)$. 

Let $X_{G}:=\beta G/G$ be the set of co-sets with respect to the right action 
of $G$.
\begin{Proposition}\label{proper}
(a)\ $E(M,\beta G)$ is a compact Hausdorff space with 
$dim\ E(M,\beta G)=dim\ M$.\\
(b)\ For every $x\in X_{G}$ there is a continuous one-to-one map
$i_{x}:M_{G}\longrightarrow E(M,\beta G)$ such that 
$\displaystyle i_{x}(M_{G})\bigcap i_{y}(M_{G})=\emptyset$ for $x\neq y$ and 
$\displaystyle E(M,\beta G)=\bigsqcup_{x\in X_{G}}i_{x}(M_{G})$.\\
(c)\ Let $f:V\longrightarrow E(M,\beta G)$ be a continuous mapping of a 
connected and locally connected space $V$. Then there
is $x\in X_{G}$ such that $f(V)\subset i_{x}(M_{G})$.\\
(d)\ For any different points $x,y\in E(M,\beta G)$ there is a connected
one-dimensional compact set $l_{x,y}\subset E(M,\beta G)$ containing them.
This set can be chosen locally connected if and only if 
$x,y\in i_{v}(M_{G})$ for some $v\in X_{G}$.
\end{Proposition}
We will also show that $i_{x}$ for $x\in X_{G}$ may turn out not to be a
homeomorphism from $M_{G}$ onto its image. 
Further, for a given structure defined on $M_{G}$ there is a related similar
``structure'' on $E(M,\beta G)$. Then we prove in these terms
an analog of the Stone-Weierstrass approximation theorem on approximation
of functions from  $C_{d}(M_{G})$ by ``Lipschitz'', ``$C^{\infty}$'' 
or ``real-analytic'' functions on $E(M,\beta G)$.
These results generalize a well-known Bernstein theorem on the uniform 
approximation of bounded uniformly continuous functions on $\Re^{n}$ by 
the entire functions bounded on $\Re^{n}$.

{\bf 2.2.}  Using the constructions of the preceding section 
we identify some naturally defined subsets of the analytical part $M_{a}$.

Let $X$ be a complex Riemann surface which admits a non-constant bounded
holomorphic function. Let $r:\Di\longrightarrow X$ be the universal covering 
over $X$. Then $\Di$ is a principle fibre bundle over $X$ with the fibre
$H:=\pi_{1}(X)$. Similarly to our main construction we determine the
fibre bundle $E(X,\beta H)$ over $X$ with the fibre $\beta H$ so that
there is a natural embedding $\Di\hookrightarrow E(X,\beta H)$ as an open
everywhere dense subset.
\begin{Th}\label{subsets}
There is an analytic injective mapping $I: E(X,\beta H)\longrightarrow M_{a}$,
i.e., $I$ is continuous and $f\circ I\circ i_{x}\in H^{\infty}$ for any 
$x\in X_{H}=\beta H/H$, $f\in H^{\infty}$.\\
For any open set $U\subset E(X,\beta H)$ the set
$I(U)\subset M_{a}$ is open. For any $x\in X_{H}$ the set
$(I\circ i_{x})(\Di)$ coincides with a Gleason part.
\end{Th}
We now state a curious corollary of the above theorem.

Let $G$ be a free group with finite or countable number of generators.
For each $\xi\in\beta G$ denote $P(\xi):=\{\xi g,\ g\in G\}$. 
\begin{C}\label{free}
Assume that 
$\overline{P(\xi_{1})}\cap\overline{P(\xi_{2})}\neq\emptyset$
for $\xi_{1},\xi_{2}\in\beta G$. Then one
of these sets contains the other.
\end{C}
\begin{R}\label{dya}
{\rm Our proof of Corollary \ref{free} is analytic.
It would be interesting to prove this result by combinatorial methods and 
to identify a class of discrete groups satisfying the above property. }
\end{R}

Let now $Y\subset Z$ be a complex Riemann surface compactly embedded into a 
Riemann surface $Z$. Assume that $H:=\pi_{1}(Y)\cong\pi_{1}(Z)$ and 
the boundary of $Y$ is a union of continuous curves. We introduce 
an analog of the disk algebra.
Let $r_{Z}:\Di\longrightarrow Z$ and $r_{Y}:\Di\longrightarrow Y$ be universal
coverings over $Z$ and over $Y$. According to the Covering Homotopy Theorem and
our assumption, there is an analytic embedding $a:\Di\longrightarrow\Di$
equivariant with respect to actions of $\pi_{1}(Y)$ and $\pi_{1}(Z)$,
i.e., $a(gz)=g(a(z))$, $z\in\Di$, $g\in\pi_{1}(Y)=\pi_{1}(Z)$.
For each open $U\subset Z$ satisfying $\overline{Y}\subset U$
consider the algebra $H_{U}$ of bounded analytic functions on 
$r_{Z}^{-1}(U)$.
Then we define $A_{Y}\subset H^{\infty}$ as closure in $C(\Di)$ 
of the algebra generated by all $a^{*}(H_{U})$.
\begin{Th}\label{algeb}
The maximal ideal space $M(A_{Y})$ of $A_{Y}$ is isomorphic to 
$E(\overline Y,\beta H)$.\\
Furthermore, $dim\ M(A_{Y})=2$ and the problem of the complementation of 
matrices formulated for algebra $A_{Y}$ always 
admits a positive solution.
\end{Th}
%==========================================
\sect{\hspace*{-1em}. Proof of Theorem \ref{new}.}
We will prove Theorem \ref{new} in a more general setting. 

Let $b$ be a $(k\times n)$-matrix with entries in a commutative 
Banach algebra $B$ of complex-valued functions defined on the maximal ideal 
space $M(B)$. 
Assume  that condition (\ref{matr}) holds at each point of $M(B)$ for the 
family of minors of $b$ of order $k$, that is, these
minors do not belong together to any maximal ideal. Assume also that
$dim\ M(B)\leq 2$. 
\begin{Proposition}\label{dimen}
There is an invertible $(n\times n)$-matrix
with entries in $B$ which completes $b$.
\end{Proposition}
According to Theorem 3 of [L], it suffices to complete $b$ in the
class of continuous matrix-functions on $M(B)$. Thus we have to
find an $(n\times n)$-matrix
$\widetilde b=(\widetilde b_{ij})$ with entries
$\widetilde b_{ij}$ in $C(M(B))$ so that $det\ \widetilde b=1$ and 
$\widetilde b_{ij}=b_{ij}$ for $1\leq i\leq k$. 
The matrix $b$ determines a trivial subbundle $\xi$ of complex rank 
$k$ in the trivial vector bundle $\theta^{n}=M(B)\times\Co^{n}$. Let 
$\eta$ be an additional to $\xi$ subbundle of $\theta^{n}$. Then clearly
$\xi\oplus\eta$ is topologically trivial. We will prove that
$\eta$ is also topologically trivial. Then a trivialization
$s_{1},s_{2},...,s_{n-k}$ of $\eta$ determines the required complement of $b$.

We begin with the following 
\begin{Lm}\label{triv}
Let $X$ be a Hausdorff compact of dimension $\leq 2$ and $p:E\longrightarrow X$
be a locally trivial continuous vector bundle over $X$ of complex rank $m$.
If $m\geq 2$ then there is a continuous nowhere vanishing section of $E$.
\end{Lm}
{\bf Proof.}
According to the general theory of vector bundles, 
see e.g., [Hu] Ch. 3, Th. 5.5,
there is a continuous map $f$ of $X$ into a Stiefel manifold
$V(m,k)$ of $m$-dimensional subspaces of $\Co^{k}$ such that for the 
canonical vector bundle
$\gamma_{m}$ over $V(m,k)$ one has $f^{*}\gamma_{m}=E$. Recall that
$\gamma_{m}$ is defined as follows: the fibre over an $x\in V(m,k)$ is
the $m$-dimensional complex space representing $x$.
We think of $V(m,k)$ as a compact submanifold of some $\Re^{s}$. Then $f$ is 
defined by a finite family of continuous on $X$ functions, 
$f=(f_{1},...,f_{s})$.
For $X'=f(X)$ let $U'$ be a compact  polyhedron containing $X'$ such 
that $\gamma_{m}$ has a continuous extension to a vector bundle over $U'$. 
Let us consider the inverse limiting system determined by all possible
finite collections of continuous on $X$ functions such that the last $s$
coordinates of any such collection are coordinates of the mapping $f$
(in the same order). For any such family $f_{\alpha}:=
(f_{1},...,f_{\alpha-s},f)$ denote 
$X_{\alpha}=f_{\alpha}(X)\subset\Re^{\alpha}$. Further, let 
$p_{\alpha}^{\beta}:X_{\beta}\longrightarrow X_{\alpha}$,
$\beta\geq\alpha$, be the mapping induced by the natural projection
$\Re^{\beta}\longrightarrow\Re^{\alpha}$ onto the first 
$\alpha$ coordinates. Let $U_{\alpha}$ be a compact polyhedron containing
$X_{\alpha}$ defined in a small open neighbourhood of $X_{\alpha}$ so that
the inverse limit of the system $(U_{\alpha},p_{\alpha}^{\beta})$ coincides
with $X$ and  the mapping $p_{0}^{\alpha}:\Re^{\alpha}\longrightarrow
\Re^{s}$ maps $U_{\alpha}$ to $U'$. Consider the bundle
$E_{\alpha}:=(p_{0}^{\alpha})^{*}(\gamma_{m})$ over $U_{\alpha}$. Then 
clearly, $(p_{\alpha}^{\beta})^{*}(E_{\alpha})=E_{\beta}$, $\beta\geq\alpha$,
and the pullback of each $E_{\alpha}$ to $X$ coincides with $E$. Now the
Euler class $e_{\alpha}$ of each $E_{\alpha}$ is an element of the
\v{C}ech cohomology group
$H^{2m}(U_{\alpha},\Z)$.
This class equals 0 if and only if $E_{\alpha}$ has a nowhere vanishing
section. By the fundamental property of characteristic classes, see e.g.
[Hi], $(p_{\alpha}^{\beta})^{*}(e_{\alpha})=e_{\beta}$, $\beta\geq\alpha$.
Further, it is well known, see e.g. [Br], Ch. II, Corol. 14.6, that 
$$
\lim_{\longrightarrow}(p_{\alpha}^{\beta})^{*}H^{k}(U_{\alpha},\Z)=
H^{k}(X,\Z),\ \ \ k\geq 0.
$$
Since $dim\ X\leq 2$, Theorem 37-7 and Corollary 36-15 of [N] imply
$H^{k}(X,\Z)=0$ for $k>2$. But real rank of $E$ is $2m\geq 4$ therefore
the Euler class $\displaystyle 
e=\lim_{\longrightarrow}(p_{\alpha}^{\beta})^{*}(e_{\alpha})$
of $E$ equals 0. From here and the above formula it follows
that there is some $\beta$ such that $e_{\beta}=0$. In particular,
$E_{\beta}$ has a continuous nowhere vanishing section $s$. Then its
pullback to $X$ determines the required section of $E$.

The lemma is proved.\ \ \ \ \ $\Box$

From the lemma it follows immediately that $E$ is isomorphic to
the direct sum $E_{m-1}\oplus E'$, where $E_{m-1}=X\times\Co^{m-1}$ and
$E'$ is a vector bundle over $X$ of complex rank 1. In particular, if
the first Chern class $c_{1}(E)\in H^{2}(X,\Z)$ of $E$ equals 0, then 
$E'=X\times\Co$ and $E=X\times\Co^{m}$.
In our case $dim\ M(B)\leq 2$, $\theta^{n}=\xi\oplus\eta$ and $\xi$ is
topologically trivial. Then for the first Chern class of $\theta^{n}$ we
have
$$
0=c_{1}(\theta_{n})=c_{1}(\xi)+c_{1}(\eta)=c_{1}(\eta).
$$
Therefore the above lemma implies that $\eta$ is a 
topologically trivial bundle.

The proposition is proved.\ \ \ \ \ $\Box$

We proceed with the proof of Theorem \ref{new}.
Let $M(A)$ be the maximal ideal space of $A$. Since $\nu(H^{\infty})$ is dense
in $A$, the transpose map
$$
\nu^{*}:M(A)\longrightarrow M(H^{\infty})\ \ (\nu^{*}(\phi)=\phi\circ\nu)
$$
establishes a homeomorphism from $M(A)$ onto some closed set
$E\subset M(H^{\infty})$. According to the result of Su'{a}rez,
$dim\ M(H^{\infty})=2$ which implies $dim\ M(A)=dim\ E\leq 2$. Further,
matrix $a$ can be thought of as a matrix-function defined on $M(A)$ whose
minors of order $k$ satisfy inequality (\ref{matr}) at each point of $M(A)$.
Thus the required statement follows from Proposition \ref{dimen} .  

The theorem is proved.\ \ \ \ \ $\Box$
%================================================
\sect{\hspace*{-1em}. Proof of Theorem \ref{te1}.}
Let $B\subset\Re^{n}$ be a Euclidean ball and $X$ be a discrete countable
set. For points $v,w\in B\times X$, $v=(v_{1},v_{2})$, $w=(w_{1},w_{2})$ 
we define the semi-metric $r(v,w):=|v_{1}-w_{1}|$, where
$|.-.|$ is the Euclidean distance on $B$.
Let $f$ be a bounded function on $B\times X$ uniformly continuous
with respect to $r$. 
\begin{Lm}\label{main}
There is a continuous function $\widehat f$ on $B\times\beta X$  such
that $\widehat f|_{B\times X}=f$ and $\sup_{B\times\beta X}|\widehat f|=
\sup_{B\times X}|f|$.
\end{Lm}
{\bf Proof.}
For any $v\in B$ function $f_{v}(x):=f(v,x)$ belongs to $l^{\infty}(X)$.
By definition, there is a continuous function $\widehat f_{v}$ on $\beta X$
which extends $f_{v}$. We set 
$\widehat f(v,\xi):=\widehat f_{v}(\xi)$ and prove that
$\widehat f$ is continuous. In fact, let us take a point 
$w=(v,\xi)\in B\times\beta X$ and a number $\epsilon>0$.
By definition of $f$, there is $\delta>0$ such that for any pair of points
$v_{1}=(v,x)$ and $v_{2}=(v',x)$ from $B\times X$ with
$|v-v'|<\delta$ one has $|f(v_{1},x)-f(v_{2},x)|<\epsilon/3$. Define an open
neighbourhood $U_{v}$ of $v\in B$ by $U_{v}:=\{v'\in B\ : |v-v'|<\delta\}$.
Further, by definition of $\widehat f_{v}$, there is an open
neighbourhood $U_{\xi}\subset\beta X$ of $\xi$ such that for any $\eta\in
U_{\xi}$ we have $|\widehat f_{v}(\eta)-\widehat f_{v}(\xi)|<\epsilon/3$.
Consider $U_{w}:=U_{v}\times U_{\xi}$. Then $U_{w}$ is
an open neighbourhood of $w\in B\times\beta X$. Note that
$f_{v}-f_{v'}$ is a function from $l^{\infty}(X)$ and 
for any $v'\in U_{v}$ its supremum norm $<\epsilon/2$. This implies that
$|\widehat f_{v}(\eta)-\widehat f_{v'}(\eta)|<\epsilon/2$ for each
$\eta\in\beta X$. In particular, for any $(x,\eta)\in U_{w}$ we have
$$
|\widehat f(x,\eta)-\widehat f(v,\xi)|\leq |\widehat f_{x}(\eta)-\widehat
f_{v}(\eta)|+|\widehat f_{v}(\eta)-\widehat f_{v}(\xi)|<\epsilon \ .
$$
This inequality shows that $\widehat f$ is continuous at
$w\in B\times\beta X$.\ \ \ \ \ $\Box$

Assume now that $B\subset\Co^{n}$ is a complex Euclidean ball and $f$ is a 
bounded holomorphic function on $B\times X$. Then $\{f(\cdot,x)\}_{x\in X}$
is a uniformly bounded family of holomorphic on $B$ functions and therefore
normality implies that $f$ satisfies conditions of Lemma \ref{main}
on every $K\times X$, where $K\subset B$ is a compact,
i.e., $f$ is uniformly continuous with respect to semi-metric $r$. Let 
$\widehat f\in C(B\times\beta X)$ be extension of $f$ determined in 
Lemma \ref{main}.
\begin{Lm}\label{hol}
For any $\xi\in\beta X$ function $\widehat f|_{B\times\xi}$ is holomorphic.
\end{Lm}
{\bf Proof.} Let $O\subset B$ be an open subset with $\overline O\subset B$.
It suffices to prove that $\widehat f|_{O\times\xi}$ is holomorhic for any $O$.
Since $\widehat f$ is uniformly continuous on the compact
$\overline O\times\beta X$, for any $\epsilon>0$ there is 
$x_{\epsilon}\in X$ such that 
$\sup_{z\in\overline O}|\widehat f(z,\xi)-f(z,x_{\epsilon})|<\epsilon$. 
In particular,
$\widehat f(\cdot,\xi)|_{O}$ is limit in $C(O)$ of the sequence 
$\{f(\cdot,x_{1/k})\}_{k\geq 1}$
of bounded holomorphic on $O$ functions. Thus $\widehat f|_{O\times\xi}$ 
is also holomorphic.

The lemma is proved.\ \ \ \ \ $\Box$
\begin{R}\label{lip}
{\rm Similar arguments show that if $f$ is a bounded on $B\times X$
function, Lipschitz with respect to $r$ then 
$\widehat f|_{B\times\xi}$ is Lipschitz for any $\xi\in\beta X$. 
However, a similar statement for a uniformly bounded on $B$ 
family of real analytic or smooth functions is false.}
\end{R}
{\bf Proof of Theorem \ref{te1}}.\\
Let $p:M_{G}\longrightarrow M$ be a regular covering over a compact
manifold $M$ and $f\in C_{d}(M_{G})$. 
Recall that, by constructions of $E(M,\beta G)$ and $M_{G}$
(see Section 2), for $\widetilde V=\sqcup_{i\in I}U_{i}\times \beta G$
there is a continuous surjective map $\widetilde s:\widetilde V
\longrightarrow E(M,\beta G)$ such that for  
$V=\sqcup_{i\in I}U_{i}\times G\subset\widetilde V$ restriction
$s:=\widetilde s|_{V}$ is a continuous map onto $M_{G}$.
Further, recall that $U_{i}$ is homeomorphic to a Euclidean ball $B$ and 
this homeomorphism determines a natural homeomorphism of 
$U_{i}\times G$ and $B\times G$. Denote by $r_{i}$ a semi-metric on
$U_{i}\times G$ which is the pullback of $r$ by the above homeomorphism.
Replacing, if necessary, in our construction the covering 
${\cal U}=\{U_{i}\}_{i\in I}$ by its finite subcovering and taking into 
account that metric $d$ is invariant with respect to the action of 
$G$ on $M_{G}$ and compatible with topology of $M_{G}$, we may assume 
without loss of generality that $(f\circ s)|_{U_{i}\times G}$ is bounded
uniformly continuous with respect to $r_{i}$.
Now according to Lemma \ref{main}, for every 
$f_{i}:=(f\circ s)|_{U_{i}\times G}$ there
is its extension $\widehat f_{i}\in C(U_{i}\times\beta G)$.
Further, the equivalence relation on $V$ implies 
$$
f_{i}(z,g)=f_{j}(z,gg_{ij})\ \ \ \ \ (z\in U_{i}\cap U_{j},\ g\in G).
$$
Then by the definition of $\widehat f_{i}$, $i\in I$,
for a fixed $z\in U_{i}\cap U_{j}$  
$$
\widehat f_{i}(z,\xi)=\widehat f_{j}(z,\xi g_{ij})\ \ \ \ \
(\xi\in\beta G) .
$$
This shows that there is a continuous function $\widehat f$ on
$E(M,\beta G)$ so that $(\widehat f\circ\widetilde s)|_{U_{i}\times\beta G}=
\widehat f_{i}$ for any $i\in I$. Moreover, $\widehat f|_{M_{G}}=f$.
Since $M_{G}$ is an open everywhere dense subset of $E(M,\beta G)$,
we also have
$$
\sup_{E(M,\beta G)}|\widehat f|=\sup_{M_{G}}|f|\ .
$$
Thus we obtain a continuous extension of $C_{d}(M_{G})$.
It remains to prove that the extended algebra coincides with $C(E(M,\beta G))$.

Note, first, that since the construction of $E(M,\beta G)$ does not depend on 
the choice of covering ${\cal U}$ of $M$, we can consider that all
$U_{i}\subset M$ are compacts. Consequently,
$E(M,\beta G)$ is a compact as continuous image of the compact 
$\widetilde V$. By definition, the base of topology on $E(M,\beta G)$
consists of the sets $\widetilde s(V\times H)$, where 
$V$ is an open subset of some $U_{i}$ and 
$H$ is a clopen (closed and open) subset of $\beta G$, i.e., closure in 
$\beta G$ of a subset of $G$. 

Let $h$ be a continuous function on $E(M,\beta G)$. Then it is bounded and
uniformly continuous. We have to prove that $h|_{M_{G}}\in C_{d}(M_{G})$, 
i.e., $h$ is bounded uniformly continuous with respect to the metric $d$. 
Let us fix $\epsilon>0$. We will point out a $\delta=\delta(\epsilon)>0$ such 
that for all $x,y$ satisfying $d(x,y)<\delta$ one has  $|h(x)-h(y)|<\epsilon$.
\\
From uniform continuity of $h$ and compactness of $E(M,\beta G)$ it follows
that for each $x\in M$ there is a finite open covering of the 
fibre $\widetilde p^{-1}(x)\subset E(M,\beta G)$ by sets $Y_{i,x}$ such
that\\
(a)\ $|h(v)-h(w)|<\epsilon$ \ \ \ \ \ $(v,w\in Y_{i,x})$;\\
(b)\  $Y_{i,x}=\widetilde s(W_{x}\times \overline{H_{i}})$,
where $W_{x}\subset M$ is an open neighbourhood of $x$ homeomorphic to a ball
and $H_{i}\subset G$.\\
In particular, $\cup_{i}Y_{i,x}=\widetilde p^{-1}(W_{x})$ and 
$Y_{x}=(\cup_{i}Y_{i,x})\cap M_{G}$ coincides with $p^{-1}(W_{x})$ and
if $Z\subset Y_{x}$ is connected and $v,w\in Z$ then
$|h(v)-h(w)|<\epsilon$. Let $\{Y_{x_{j}}\}_{j\in J}$ be a finite covering
of $E(M,\beta G)$ and $\{W_{x_{j}}\}_{j\in J}$ be the corresponding covering 
of $M$. 
For any connected component $Z$ of $Y_{x_{j}}$ denote by $a_{j}$ its
diameter with respect to $d$ and set $a:=\max_{j\in J}a_{j}$. Clearly
$0<a<\infty$. Since $G$ acts discretely on $M_{G}$, $d$ is invariant 
with respect to this action and compatible with topology on $M_{G}$, and the 
quotient $M=M_{G}/G$ is a compact, we obtain 
$$
\inf_{g\in G,g\neq e}d(x,gx)=b>0.
$$
Doing, if necessary, a refinement of the covering $\{Y_{x_{j}}\}$ we can 
assume that $b-a>0$.
Further, $d^{*}(p(x),p(y)):=\inf_{g\in G}d(x,gy)$ is a metric
on the compact $M$ compatible with its topology. Then there is a positive
number $\delta=\delta(\epsilon)<b-a$ such that every ball 
$B_{\delta}\subset M$ of radius $\delta$ with respect to $d^{*}$ belongs to 
one of $W_{x_{j}}$. We state that 

for any $x,y\in M_{G}$ with $d(x,y)<\delta$ the inequality
$|h(x)-h(y)|<\epsilon$ holds.\\
In fact, according to our construction $d^{*}(p(x),p(y))<\delta$ and so there 
is $W_{x_{j}}$ containing $p(x)$ and $p(y)$. Let $Z$ be a connected component 
of $p^{-1}(W_{x_{j}})$ containing $x$. Then it also contains $y$. 
For otherwise, let $gy\in Z$ with $g\neq e$. Then 
$$
\delta>d(x,y)\geq d(y,gy)-d(gy,x)\geq b-a\ .
$$
This contradicts to the choice of $\delta$ and proves that $g=e$ and
so $|h(x)-h(y)|<\epsilon$ by the choice of $\{Y_{x_{j}}\}_{j\in J}$.

The above arguments show that $h|_{M_{G}}$ is a bounded uniformly continuous 
with respect to $d$ function.

Thus we proved that algebra $C_{d}(M_{G})$ admits a continuous extension 
(preserving the norm) to the compact $E(M,\beta G)$ and the extended algebra 
coincides with $C(E(M,\beta G))$. Therefore the maximal ideal space of 
$C_{d}(M_{G})$ is $E(M,\beta G)$.

The proof of the theorem is complete. \ \ \ \ \ $\Box$
%======================================================
\sect{\hspace*{-1em}. Properties of the Space $E(M,\beta G)$.}
{\bf 5.1.} All results of this section, except for compactness of
$E(M,\beta G)$, are valid for a non-compact $M$. We recall
\begin{D}\label{dimensi}
For a normal space $X$, $dim\ X\leq n$ if every finite
open covering of $X$ can be refined by an open covering whose order $\leq
n+1$. If $dim\ X\leq n$ and the statement $dim\ X\leq n-1$ is false,
we say $dim\ X=n$. For the empty set, $dim\ \emptyset=-1$.
\end{D}
{\bf Proof of Proposition \ref{proper}}.
{\bf (a)}
We already proved that $E(M,\beta G)$ is a compact Hausdorff space. By 
definition, it admits a finite covering by subsets homeomorphic to
$B\times\beta G$, where $B$ is a closed Euclidean ball in $\Re^{n}$ and
$n=dim\ M$. We will show that $dim\ B\times\beta G=n$. Then $E(M,\beta G)$
is a finite union of compact sets of dimension $n$ and therefore
$dim\ E(M,\beta G)=n$, see, e.g., [N], Th. 9-10.
 
Since $B\times \{e\}\subset B\times\beta G$, we obviously have 
$dim\ B\times\beta G\geq dim\ B=n$. Further, Theorem 26-4 of [N] implies 
$$
dim\ B\times\beta G\leq dim\ B+dim\ \beta G.
$$
Finally by the theorem of Vendenisov, see, e.g., [N], Th. 9-5, $dim\ \beta G=
dim\ G=0$.  Combining the above inequalities we obtain
$dim\ B\times\beta G=n$. This proves part (a).\\
{\bf (b)}
We begin with the following
\begin{Lm}\label{fixed}
For any $\xi\in\beta G$ and $g\in G$, $g\neq e$, one has
$\xi g\neq\xi$ and $g\xi\neq\xi$.
\end{Lm}
{\bf Proof.}
We will prove this result for the right action of $G$ on $\beta G$. 
For the left action the proof is similar.
Assume, to the contrary, that $\xi h=\xi$ for some $h\in G,\ h\neq e$, 
$\xi\in\beta G$. Set $H:=\{h^{n},\ n\in\Z\}\subset G$. Consider right
co-sets $\{G_{\alpha}:=g_{\alpha}H\}_{\alpha\in A}$ by the action of $H$. 
Each $G_{\alpha}$ is disjoint union of $G_{\alpha,0}$ and 
$G_{\alpha,1}$, where $G_{\alpha,0}:=\{g_{\alpha}h^{n},\ n\in 2\Z\}$ and
$G_{\alpha,1}:=\{g_{\alpha}h^{n},\ n\in 2\Z+1\}$. Let 
$G_{0}=\cup_{\alpha}G_{\alpha,0}$ and $G_{1}=\cup_{\alpha}G_{\alpha,1}$.
Then $G_{0}\cap G_{1}=\emptyset$ and $G_{0}\cup G_{1}=G$. Moreover,
$G_{0}h=G_{1}$ and $G_{1}h=G_{0}$. Obviously, $\beta G$ is disjoint union
of closures $\overline{G_{0}}$ and $\overline{G_{1}}$. Assume, e.g., that
$\xi\in\overline{G_{0}}$. Then $\xi h\in\overline{G_{1}}$ and we obtain a
contradiction.

The lemma is proved.\ \ \ \ \ $\Box$
 
Let now $x\in X_{G}:=\beta G/G$ and $\xi\in\beta G$ be a point representing
$x$. Further, let $\widetilde V:=\sqcup_{i} U_{i}\times\beta G$ be the same
as in the construction of $E(M,\beta G)$ (see Section 2). Consider set
$\widetilde V_{x}:=\sqcup_{i}U_{i}\times\xi G\subset\widetilde V$. For
$z\in U_{i}\cap U_{j}$ we identify points $z\times\xi g\in U_{i}\times
\xi G$ and $z\times\xi gg_{ij}\in U_{j}\times\xi G$.  Here 
$g=\{g_{ij}\}\in Z^{1}({\cal U}, G)$ is a cocycle determining
$E(M,\beta G)$. Now according to Lemma \ref{fixed} there is a natural 
continuous one-to-one map of
$V:=\sqcup_{i} U_{i}\times G$ into $\widetilde V$ sending
$z\times g\in U_{i}\times G$ to $z\times\xi g\in U_{i}\times\xi G$. This
map commutes with the corresponding equivalence relations 
determining $M_{G}$ and $E(M,\beta G)$, respectively. Therefore it determines
a  continuous one-to-one map
$i_{x}:M_{G}\longrightarrow E(M,\beta G)$. Clearly 
$i_{x}(M_{G})\cap i_{y}(M_{G})=\emptyset$ for $x\neq y$ and 
$E(M,\beta G)=\sqcup_{x\in X_{G}}i_{x}(M_{G})$. It is also worth noting that
locally $i_{x}$ is an embedding. This proves part (b).\\
{\bf (c)}
Let $f:V\longrightarrow E(M,\beta G)$ be a continuous mapping of a
connected and locally connected space $V$. We prove that there is
$x\in X_{G}$ such that $f(V)\subset i_{x}(M_{G})$. 
For every $v\in V$ consider $f(v)\in E(M,\beta G)$. Let $U_{f(v)}$ be an
open neighbourhood of $f(v)$ homeomorphic to $O\times H$, where
$O\subset M$ and $H\subset\beta G$ are open and $U_{f(v)}\subset\widetilde p
^{-1}(O)$ (see description of topology on $E(M,\beta G)$). Then  there is a 
connected open neighbourhood $U_{v}\subset V$ of $v$ such that 
$f(U_{v})\subset U_{f(v)}$. Let $\pi_{O}$ be the projection from 
$\widetilde p^{-1}(O)$ onto a fibre of $E(M,\beta G)$ inside
(homeomorphic to $\beta G$). Then 
$\pi_{O}\circ f$ maps $U_{v}$ onto a subset of $\beta G$. Since 
the above mapping is continuous, $U_{v}$ is connected and $\beta G$ is 
totally disconnected ($dim\ \beta G=0$), we have 
$(\pi_{O}\circ f)|_{U_{v}}$ is a constant map.
In particular, there is $x\in X_{G}$ such that $f(U_{v})\subset i_{x}(M_{G})$.
Further, since $V$ is a connected space, there is $U_{v'}$ with similar
properties such that $U_{v'}\cap U_{v}\neq\emptyset$. Let $x'\in X_{G}$ be
such that $f(U_{v'})\subset i_{x'}(M_{G})$. Then the above condition implies
$i_{x'}(M_{G})\cap i_{x}(M_{G})\neq\emptyset$. Therefore according to 
Proposition \ref{proper} (b), $x'=x$. Continuing by induction we finally 
obtain $f(V)\subset i_{x}(M_{G})$.\ The proof of part (c) is complete.\\
{\bf (d)}
Let $x,y\in E(M,\beta G)$ be different points and $\widetilde p(x),
\widetilde p(y)$ their images in $M$. Let $\gamma_{x,y}$ be a path connecting 
$\widetilde p(x)$ with $\widetilde p(y)$ and 
$\{\gamma_{i}\}_{i\in I}$ be a family of closed paths with origin at $y$ 
which generates fundamental group $\pi_{1}(M)$. Consider compact space 
$E_{x,y}:=\widetilde p^{-1}(\gamma_{x,y}\cup (\cup_{i\in I}\gamma_{i}))
\subset E(M,\beta G)$. Arguments similar to those used in (a) applied to
$E_{x,y}$ show that $dim\ E_{x,y}=1$.
Further, by the Covering Homotopy Theorem $E_{x,y}\cap M_{G}$ is an open 
arcwise connected everywhere dense subset of $E_{x,y}$. This implies that
$E_{x,y}$ is connected. Let $l_{x,y}$ be a compact 
connected subset of $E_{x,y}$ containing $x$ and $y$. Then it is 
one-dimensional. For otherwise, $dim\ l_{x,y}=0$ and therefore $l_{x,y}$ is 
totally disconnected that contradicts to its choice. \\
Assume, in addition, that $l_{x,y}$ is locally connected. Then according to
part (c), there is $z\in X_{G}$ such that $l_{x,y}\subset i_{z}(M_{G})$. 
In particular, $x,y\in i_{z}(M_{G})$. The converse to this statement is 
obvious: if $x,y\in i_{z}(M_{G})$ for some $z\in X_{G}$ then there is a 
compact path $\gamma\subset M_{G}$ such that $i_{z}(\gamma)$ contains $x$ and 
$y$. (Note also that $i_{z}|_{\gamma}:\gamma\longrightarrow i_{z}(\gamma)$ is
a homeomorphism).
This proves part (d).

The proposition is proved.\ \ \ \ \ $\Box$
\begin{R}\label{ends}
{\rm There is a natural continuous surjective map $e$ of 
$E(M,\beta G)\setminus M_{G}$ onto
space $E(M_{G})$ of ends of $M_{G}$ (it is a totally disconnected compact 
Hausdorff space).
If $x,y\in e^{-1}(a)$ for some $a\in E(M_{G})$ then one can also choose 
$l_{x,y}$ inside $e^{-1}(a)$.}
\end{R}
{\bf 5.2.} In this section we study maps 
$i_{x}:M_{G}\longrightarrow E(M,\beta G)$,\ $x\in X_{G}$, which are embeddings
(i.e., homeomorphisms from $M_{G}$ onto $i_{x}(M_{G})$) and are not embeddings.

Let $Y=p^{-1}(y)$ be fibre of the bundle $p:M_{G}\longrightarrow M$ over a 
point $y\in M$. Consider set $Z:=i_{x}(Y)\subset E(M,\beta G)$ and denote
$Cl(Z)$ the set of limit points of $Z$.
\begin{Proposition}\label{embed}
$i_{x}$ is an embedding if and only if $Cl(Z)\cap Z=\emptyset$.
\end{Proposition}
{\bf Proof.}
Assume that $i_{x}$ is an embedding. Since
$Y\subset M_{G}$ is discrete, its image $Z$ is discrete in 
$E(M,\beta G)$. In particular, $Cl(Z)$ does not contain points of $Z$.\\
Conversely, assume, to the contrary, that $Cl(Z)\cap Z=\emptyset$
but there is a discrete sequence $\{x_{k}\}_{k\geq 1}$ of points 
of $M_{G}$ so that the limit set $Cl(X)$ of $X:=i_{x}(\{x_{k}\}_{k\geq 1})$ 
has non-empty intersection with $i_{x}(M_{G})$. 
\begin{Lm}\label{inprop}
There is a continuous map $\widehat{i_{x}}:E(M,\beta G)\longrightarrow
E(M,\beta G)$ which extends $i_{x}$, i.e., $\widehat{i_{x}}|_{M_{G}}=i_{x}$.
\end{Lm}
{\bf Proof.}
For every continuous on $E(M,\beta G)$ function $f$ consider
the function $f':=f\circ i_{x}$ on $M_{G}$. Then it follows easily from
uniform continuity of $f$ that $f'\in C_{d}(M_{G})$, i.e., it is bounded and
uniform 
continuous with respect to the metric $d$. In particular, $f'$ admits a
continuous extension $\widehat{f'}$ to $E(M,\beta G)$. Clearly, 
the correspondence $f\mapsto\widehat{f'}$ 
is a continuous morphism of Banach algebra $C(E(M,\beta G))$ 
to itself. In turn, this morphism determines the continuous transpose map
$\widehat{i_{x}}:E(M,\beta G)\longrightarrow E(M,\beta G)$ which extends
$i_{x}$. 

The lemma is proved.\ \ \ \ \ $\Box$

Let $X'\subset E(M,\beta G)$ be closure of $\{x_{k}\}_{k\geq 1}$.
Obviously, $\widehat{i_{x}}(X')=\overline{X}$.
By our assumption, there is a limit point $s\in X'$ of $\{x_{k}\}_{k\geq 1}$
such that $\widehat{i_{x}}(s)=i_{x}(z)$ for some $z\in M_{G}$. Assume also 
that $s$ belongs to $i_{t}(M_{G})\subset E(M,\beta G)$, for some $t\in X_{G}$.
Further, by definition of $E(M,\beta G)$ closure of $Y$ in $E(M,\beta G)$ has
non-empty intersection with $i_{t}(M_{G})$. Let $v\in\overline{Y}
\cap i_{t}(M_{G})$.
Then according to Proposition \ref{proper} (c), the set
$(\widehat{i_{x}}\circ i_{t})(M_{G})$ belongs to some $i_{l}(M_{G})$. But this 
set contains also the point $i_{x}(z)$ that implies $l=x$. In particular,
$\widehat{i_{x}}(v)\in i_{x}(M_{G})$. Since, by definition, 
$\widehat{i_{x}}(v)\in Cl(Z)$, we have that  
$Cl(Z)\cap i_{x}(M_{G})\neq\emptyset$. But
closure $\overline{Z}$ lies in the fibre $\widetilde p^{-1}(y)$ over 
$y\in M$.
Thus  $Cl(Z)\cap Z\neq\emptyset$ that contradicts
to our assumption. Therefore for any discrete sequence of points of $M_{G}$
its image under $i_{x}$ is discrete in $E(M,\beta G)$. This 
proves that $i_{x}$ is an embedding.

The proposition is proved.\ \ \ \ \ $\Box$

We will show now how to construct maps
$i_{x}:M_{G}\longrightarrow E(M,\beta G)$ which are embeddings. 

Let us fix a system $\{G_{i}\}$ of finite subsets of the 
group $G$ so that $G_{i}^{-1}=G_{i}$ and $G_{i}\subset G_{i+1}$ and 
$\cup G_{i}=G\setminus\{e\}$. 
Construct by induction the following sequence $\{g_{i}\}$ of elements of $G$:

$g_{1}$ is an arbitrary element and then we take $g_{n}$ as an 
element from complement to 
$$
\left(\bigcup_{i=1}^{n}g_{1}G_{i}\right)\bigcup
\left(\bigcup_{i=1}^{n}g_{2}G_{i}\right)\bigcup...\bigcup
\left(\bigcup_{i=1}^{n}g_{n-1}G_{i}\right)\ .
$$
Then clearly for any $i,n$ we have $g_{n}\not\in g_{i}G_{i}$. Let now
$\xi\in\beta G$ be a limit point of $\{g_{i}\}$ and $x:=x(\xi)\in X_{G}$ be
the element corresponding to $\xi$. 
\begin{Proposition}\label{hom}
$i_{x}$ determines an embedding of $M_{G}$ into $E(M,\beta G)$.
\end{Proposition}
{\bf Proof.}
According to Proposition \ref{embed} it suffices to prove that 
$\{\xi G\}$ is a discrete subset of $\beta G$. 
(Here we identify fibres of $E(M,\beta G)$ with $\beta G$.) Note also
that it suffices to prove that $\xi$ itself does not belong to the limit set
of $\{\xi G\}$.
Assume, to the contrary, that there is a discrete net 
$\{h_{\alpha}\}\subset G$ such 
that $\lim_{\alpha}\xi h_{\alpha}=\xi$. Let us construct function $f$
on $G$ by 
$$
f(g_{i})=1\ \ \ i=1,2,...\ \ \ {\rm and}\ \ \ f|_{G\setminus\{g_{i}\}}\equiv 0.
$$
Let $\widehat f$ be the continuous function on $\beta G$ which
extends $f$. Then, by definition, $\widehat f(\xi)=1$. Let us prove that  
$\widehat f(\xi g)=0$ for any $g\neq e$. In fact, consider the sequence
$\{g_{i}g\}_{i\geq 1}$. Then clearly its closure contains $\xi g$. 
By definition, there
is a number $i_{0}$ such that $g\in G_{j}$ for any $j\geq i_{0}$. But the
construction of $\{g_{i}\}$ implies that $\{g_{j}g\}_{j\geq i_{0}}$ does not
intersect $\{g_{i}\}$. In particular, $f(g_{j}g)=0$ for any $j\geq i_{0}$.
Therefore we obtain $\widehat f(\xi g)=0$ for any $g\in G, g\neq e$. 
This gives the contradiction to the condition $\lim_{\alpha}
\widehat f(\xi h_{\alpha})=\widehat f(\xi)$. 

The proposition is proved.\ \ \ \ \ $\Box$

We present now non-homeomorphic maps 
$i_{x}:M_{G}\longrightarrow E(M,\beta G)$. We use the following
\begin{Lm}\label{semi}
There is a multiplication  $\cdot :\beta G\times\beta G\longrightarrow\beta G$ 
converting $\beta G$ into a semigroup. Moreover, for $g\in G$ the mapping
$\xi\mapsto \xi\cdot g$, $\xi\in\beta G$ coincides with the right action
of $G$ on $\beta G$ and 
for any $\xi\in\beta G$ the mapping 
$\eta\mapsto\eta\cdot\xi$, $\eta\in\beta G$, is continuous.
\end{Lm}
{\bf Proof.}
Let $\xi\in\beta G$. Consider the set $o(\xi):=\{\xi g\}_{g\in G}$.
The restriction of a function from $C(\beta G)$ to $o(\xi)$
determines a continuous homomorphism of $C(\beta G)$ into 
$l^{\infty}(G)$. Denote by
$f_{\xi}:\beta G\longrightarrow\beta G$ the continuous transpose map of
maximal ideal 
spaces corresponding to this homomorphism. Further, let us consider the space 
$T$ of all mappings $\beta G\longrightarrow\beta G$ with the product topology.
Then $T$ is a compact Hausdorff space.  Let $K$ be  
closure in $T$ of the group $G$ acting by right transformations on $\beta G$.
By definition, $K$ is a semigroup (in discrete topology) and the mapping 
$\eta\mapsto\eta\circ\xi$,
$\eta\in K$, is continuous for any fixed $\xi\in\beta G$. Assume now that
a net $\{g_{j}\}_{j\in J}$ of elements from $G$ converges to $k\in K$.
Then,
for any $\xi\in\beta G$ the element $k(\xi)\in\beta G$ is defined as
$$
k(\xi)=\lim_{j}\xi\cdot g_{j}.
$$
In particular, this formula applied to $\xi=e$ states that there is 
limit $\eta=\lim_{j}g_{j}\in\beta G$. But then we clearly have
$k(\xi)=f_{\xi}(\eta)$. Define now map $i:K\longrightarrow\beta G$ by
$i(k)=k(e)$. It is obviously bijective and continuous in the topology
on $K$ induced from $T$. Therefore this map is homeomorphic. Then
for $\eta,\xi\in\beta G$ we determine 
$$
\eta\cdot\xi:=i(i^{-1}(\eta)\circ i^{-1}(\xi)).
$$
The required properties of the multiplication follow easily from this
definition.

The lemma is proved.\ \ \ \ \ $\Box$

Let now $\xi\in \beta G\setminus G$ be an idempotent with respect to the above 
multiplication, i.e., $\xi\cdot\xi=\xi$ (its existence for compact
semigroups is proved, e.g., in
[F], p.159, Lemma 8.4). Then $i_{\xi}:M_{G}\longrightarrow
E(M,\beta G)$ is not embedding. In fact, the limit set of 
$\{\xi g\}_{g\in G}$ contains $\xi$. Thus the required statement 
follows from Proposition \ref{embed}.
%=====================================
\sect{\hspace*{-1em}. Approximation Theorems.}
In this section we prove a theorem on approximation of continuous functions
on $E(M,\beta G)$. Assume that the manifold
$M$ has one of the structures: Lipschitz, $C^{\infty}$, real analytic or
holomorphic. Denote by ${\cal E}$ one of the above categories.
\begin{D}\label{stru}
A complex-valued function $f\in C(E(M,\beta G))$ belongs to ${\cal E}$
if the function $f\circ i_{x}$ belongs to 
${\cal E}$ on $M_{G}$ for each $x\in X_{G}$.
\end{D}
Similarly to this definition one can define sheaves of Lipschitz, 
$C^{\infty}$, real analytic or holomorphic functions on $E(M,\beta G)$. 

Assume now that $M$ is a compact real analytic manifold and
$C_{\Re}(E(M,\beta G))$
is the Banach space of real-valued continuous functions on $E(M,\beta G)$.
\begin{Th}\label{approx}
The space of real analytic on $E(M,\beta G)$ functions  is dense in 
$\penalty -10000$ $C_{\Re}(E(M,\beta G))$.
\end{Th}
{\bf Proof.}
According to Grauert's theorem, see [Gr], 
sect.3, $M$ admits a real analytic embedding into some $\Re^{n}$; so
without loss of generality, assume that $M\subset\Re^{n}$. Further, since
any compact of $\Re^{n}$ is polynomially convex in $\Co^{n}$ with respect
to holomorphic polynomials (a consequence of the Weierstrass approximation
theorem) and $M$ is smooth, there is an open connected neighbourhood 
$U\subset\Co^{n}$ of $M$ such that:\\
(a)\ \ \ $M$ is a deformation retract of $U$;\\
(b)\ \ \ $U$ is a Stein manifold;\\
(c)\ \ \ $U$ is invariant with respect to the involution $s$ on 
$\Co^{n}$ defined by complex conjugation of coordinates.\\
According to $(a)$,
we have an isomorphism $\pi_{1}(U)\cong\pi_{1}(M)$. Then by the
Covering Homotopy Theorem there is an analytic embedding 
$i_{G}:M_{G}\longrightarrow U_{G}$, where $U_{G}$ is a regular covering
over $U$ with the covering group $G$. (Denote the covering
projection of $U_{G}$ onto $U$ by $p$ (as for $M_{G}$).)
Moreover, $U_{G}$ is also a Stein manifold. 
Clearly $i_{G}$ determines
an embedding $\widetilde i_{G}:E(M,\beta G)\longrightarrow E(U,\beta G)$ such
that for every holomorphic on $E(U,\beta G)$ function $f$ the function
$f\circ\widetilde i_{G}$ is complex-valued analytic
on $E(M,\beta G)$. (We will say that 
$\widetilde i_{G}$ is an analytic embedding.) 
Without loss of generality we may now assume that
$E(M,\beta G)\subset E(U,\beta G)$.\\
We prove that\\
(1)\ {\em for any different points} $x_{1},x_{2}\in E(M,\beta G)$ 
{\em there is a holomorphic function on}  $E(U,\beta G)$ 
{\em separating them.}\\
(2)\ {\em restrictions of holomorphic on} $E(U,\beta G)$ {\em functions 
to} $E(M,\beta G)$ {\em generate a self-adjoint with respect to complex 
conjugation algebra $A$.}\\
Then according to the classical Stone-Weierstrass theorem
$A$ is dense in complexification of $C_{\Re}(E(M,\beta G))$.
To complete the proof
of the theorem it will be sufficient to take real parts of elements of $A$,

{\bf Proof of (1)}.
Let $\widetilde p(x_{1}),\widetilde p(x_{2})\in M$ be projections of $x_{1}$ 
and $x_{2}$. Assume that $\widetilde p(x_{1})\neq\widetilde p(x_{2})$.
(if these points coincide the proof is similar).
Let $E_{x_{1}}, E_{x_{2}}\subset E(M,\beta G)$ be fibers over 
$\widetilde p(x_{1}),\widetilde p(x_{2})$. We prove that for any complex
continuous functions $f_{1},f_{2}$ determined
on $E_{x_{1}}$ and $E_{x_{2}}$, respectively, there is a holomorphic function
$f$ on $E(U,\beta G)$ such that $f|_{E_{x_{i}}}=f_{i}$, $i=1,2$.
This will imply (1).

The natural right action of $G$ on itself
determines a natural isometric action on $l^{\infty}(G)$. Let 
$g=\{g_{ij}\}\in Z^{1}({\cal U},G)$ be a cocycle determining the bundle
$p:M_{G}\longrightarrow M$ (see Section 2). Then we can also think of $g$ as a 
family of linear endomorphisms of $l^{\infty}(G)$. Using this cocycle 
$g=\{g_{ij}\}\in Z^{1}({\cal U}, GL(l^{\infty}(G)))$
we construct a Banach holomorphic vector bundle $F$ over $M$ with fibre 
$l^{\infty}(G)$  by a construction analogous to the construction of $M_{G}$
(see Section 2.1 and definitions in [B], sect.4). 
Namely, $F$ is the quotient space of disjoint union 
$\sqcup_{i\in I}U_{i}\times l^{\infty}(G)$ by the equivalence relation:
$U_{i}\times l^{\infty}(G)\ni (x,v)\sim (x,vg_{ij})\in U_{j}\times 
l^{\infty}(G)$. Here we denote $v\mapsto vg$ the action of $g$ on 
$l^{\infty}(G)$ induced by right multiplication on $G$.
Let now
$r:U\longrightarrow M$ be the retraction from property (a) for the manifold
$U$. Consider the covering $r^{-1}({\cal U})=\{r^{-1}(U_{i})\}_{i\in I}$ of 
$U$, where ${\cal U}=\{U_{i}\}_{i\in I}$ is the above covering of $M$. 
Then $U_{G}$ can be defined by cocycle 
$g=\{g_{ij}\}\in Z^{1}(r^{-1}({\cal U}), G)$ with the same elements $g_{ij}$.
Similarly to the construction of $F$ we construct
a Banach holomorphic vector bunlde $E$ over $U$ with fibre $l^{\infty}(G)$
by the cocycle $g=\{g_{ij}\}\in Z^{1}(r^{-1}({\cal U}),GL(l^{\infty}(G)))$. 
Obviously we have, $E|_{M}=F$.

Consider now a vector space $S$ of holomorphic functions $f$ on $U_{G}$ such
that $\penalty -10000$ $\sup_{z\in p^{-1}(O)}|f(z)|<\infty$ for any compact
$O\subset U$.
Then there is a natural isomorphism between $S$ and the space of
holomorphic sections of $E$. Indeed, if $f\in S$ then it can be written in
local coordinates $(z,g)$,\ $z\in r^{-1}(U_{i})\subset U,\ g\in G,\ i\in I$ as 
$f_{i}=\{f(gs_{i}(z))\}_{g\in G}$, where
$s_{i}:r^{-1}(U_{i})\longrightarrow U_{G}$ 
is a local holomorphic section of the bundle $p:U_{G}\longrightarrow U$ over 
$r^{-1}(U_{i})$. From the normality of family $f_{i}$ it follows that 
$\{(f_{i},U_{i})\}_{i\in I}$ uniquely determines a holomorphic section $s$ of 
$E$ such that $s|_{U_{i}}=f_{i}$. Further, according to Theorem \ref{te1} and
Lemma \ref{hol} any function from $S$ admits an extension as a holomorphic
function on $E(U,\beta G)$.
Conversely, for any holomorphic on $E(U,\beta G)$ function $f$ the 
function $f|_{U_{G}}\in S$.

Let ${\cal O}_{E}$ be a Banach coherent analytic sheaf of germs of local 
holomorphic sections of $E$ and ${\cal J}$ be its subsheaf
of germs of local holomorphic sections vanishing at $\widetilde p(x_{1})$ and
$\widetilde p(x_{2})$. The natural embedding  
$\widetilde p(x_{1})\sqcup\widetilde p(x_{2})\subset U$ induces surjective
homomorphism of sheaves ${\cal O}_{E}\longrightarrow {\cal O}_{x_{1},x_{2}}$, 
where ${\cal O}_{x_{1},x_{2}}$ is the sheaf of sections of $E$ over
$\widetilde p(x_{1})\sqcup\widetilde p(x_{2})$. Clearly its kernel coincides 
with ${\cal J}$ and the factor sheaf ${\cal O}_{E}/{\cal J}$ is isomorphic to
${\cal O}_{x_{1},x_{2}}$. It is well-known that the following sequence of 
cohomology of sheaves\\
$$
0\longrightarrow H^{0}(U,{\cal J})\longrightarrow H^{0}(U,{\cal O}_{E})
\longrightarrow
H^{0}(U,{\cal O}_{E}/{\cal J})\longrightarrow H^{1}(U,{\cal J})\longrightarrow
...
$$
is exact. Moreover, according to results of Bungart (see [B], sect.4),
$H^{1}(U,{\cal J})=0$. In particular, there is a holomorphic
section $f$ of $E$ whose restrictions to 
$\widetilde p(x_{1})$ and $\widetilde p(x_{2})$ coincide with 
$f_{1}|_{p^{-1}(\widetilde p(x_{1}))}$ and 
$f_{2}|_{p^{-1}(\widetilde p(x_{2}))}$, respectively.
It remains to take the holomorphic function on $U_{G}$ which represents 
$f$ and extend it to $E(U,\beta G)$.
This proves (1).

{\bf Proof of (2).}
Since $U$ is invariant with
respect to the involution $s$ defined in $(c)$ and $s$ induces identity
authomorphism of $\pi_{1}(U)$, the Covering Homotopy Theorem implies
that there is an involution $s'$ of $U_{G}$ which covers  
$s$ such that $s'|_{M_{G}}$ is trivial. 
By definition of $s'$, for each holomorphic section $f$ of $E$,  
$\overline{f\circ s'}$ is also holomorphic section of $E$ and 
its restriction to $M$ coincides with $\overline{f}|_{M}$. Identifying as in
the proof of (1) holomorphic sections of $E$ with holomorphic functions 
on $E(U,\beta G)$ we obtain that
for each holomorphic on $E(U,\beta G)$ function $f$ there is a holomorphic
on $E(U,\beta G)$ function $h$ such that 
$h|_{E(M,\beta G)}=\overline{f}|_{E(M,\beta G)}$.
This proves (2) and completes the proof of the theorem.
\ \ \ \ \ $\Box$.
\begin{R}\label{aproxim}
{\rm Similar to Theorem \ref{approx} statements are also valid for approximations by
Lipschitz or $C^{\infty}$ functions. The proof goes along 
the same lines and might be left to the reader.}
\end{R}
\begin{R}\label{bern}
{\rm Arguing as in the proof of Theorem \ref{approx} one can obtain another
proof of classical Bernstein's theorem on 
uniform approximation of bounded uniformly continuous (with respect to the
Euclidean metric) functions on $\Re^{n}$ by restrictions of entire functions 
of exponential type bounded on $\Re^{n}$. The essential component of the proof
is that $\Re^{n}$ covers a torus $\To^{n}$.}
\end{R}
%==========================================
\sect{\hspace*{-1em}. Analytical Part of Maximal Ideal Space of $H^{\infty}$.}
In this section we prove Theorem \ref{class}. Let $M$ be a complex compact 
Riemann surface of genus $g\geq 2$ and $G$ its fundamental group. Then 
$G$ acts on $\Di$ by M\"{o}bius transformations.
Consider the compact space $E(M,G)$. We begin with description
of the maximal ideal space $M(H^{\infty})$ as a quotient space of $E(M,G)$.

By definition, $\Di$ is embedded in $E(M,\beta G)$ as an open everywhere 
dense subset and every function $f\in H^{\infty}$ admits an extension to a 
holomorphic function on $E(M,\beta G)$ (see Lemma \ref{hol} and Theorem 
\ref{te1}). We identify $f$ with its extension and say that 
$H^{\infty}$ is defined on $E(M,\beta G)$.
For points $x,y\in E(M,\beta G)$ we say that $x$ is equivalent to $y$ if 
$f(x)=f(y)$ for any $f\in H^{\infty}$. Then the quotient space $E$ of 
$E(M,\beta G)$ by this equivalence relation
is a compact Hausdorff space in the factor topology defined by

$U\subset E$ is open if and only if its preimage in $E(M,\beta G)$ is open.
\\
Denote by $F:E(M,\beta G)\longrightarrow E$ the corresponding quotient map.
Then $F$ is continuous and surjective. Further, note 
that $H^{\infty}$ separates points of $\Di$ and  
if $x\in\Di$ and $y\in E\setminus\Di$ then $x$ is not equivalent to $y$. In 
fact, for coordinate function $z\in H^{\infty}$ we have $|z(x)|<1$ but 
$|z(y)|=1$.
This shows that $F$ embeds $\Di$ into $E$ as an open everywhere dense subset.
Now for every $g\in H^{\infty}$ there is a unique continuous function $g'$ 
on $E$ such that  $F^{*}g'=g$. Therefore without loss of generality 
we may assume that $H^{\infty}$ is defined on $E$. Clearly, $H^{\infty}$ 
separates points of $E$.
\begin{Proposition}\label{homeo}
$E$ is homeomorphic to the maximal ideal space $M(H^{\infty})$.
\end{Proposition}
{\bf Proof.}
Since $H^{\infty}$ is a subalgebra of $C(E)$ and for any 
$x\in E$ and $f\in H^{\infty}$ we have 
$|f(x)|\leq\sup_{\Di}|f|$, every point of $E$ determines a
continuous homomorphism of $H^{\infty}$. Therefore there is an embedding
$i:E\hookrightarrow M(H^{\infty})$. Then
$i(E)$ is a compact subset of $M(H^{\infty})$.
By the corona theorem $\Di$ is everywhere dense in 
$M(H^{\infty})$. This implies that $i$ is surjective. So $i$ is a 
homeomorphism as a continuous bijection of Hausdorff compacts.

The proposition is proved.\ \ \ \ \ $\Box$\\
{\bf Proof of Theorem \ref{class}.}
{\bf (a)}
We recall the following
\begin{D}\label{interp}
A sequence $\{z_{j}\}\subset\Di$ is an interpolating sequence if
for every bounded sequence of complex numbers $\{a_{j}\}$, there is an
$f\in H^{\infty}$ so that $f(z_{j})=a_{j}$.
\end{D}
A sequence $\{g_{j}\}\subset G$ is said to be interpolating if 
$\{g_{j}(0)\}\subset\Di$ is interpolating. 
Let $G_{in}\subset\beta G$ be the set defined by

$\xi\in G_{in}$ {\em iff there is an interpolating sequence} $
\{g_{j}\}\subset G$ {\em such that}  $\xi\in\overline{\{g_{j}\}}$.\\
In Lemma \ref{Gin} we prove that $G_{in}$ is invariant with respect to 
the right 
action of $G$. Therefore similar to the construction of  $E(M,\beta G)$ one can
construct a locally trivial bundle $E(M,G_{in})$ over $M$ with 
fibre $G_{in}$ (by the same cocycle $g\in Z^{1}({\cal U},G)$ determining 
$E(M,\beta G)$). Moreover, since $G_{in}\subset\beta G$, this bundle
has a natural embedding into $E(M,\beta G)$. We will show that
$F|_{E(M,G_{in})}$ is a homeomorphism onto the analytical part 
$M_{a}\subset M(H^{\infty})$. This will complete the proof of (a).
In our proof we use the following facts.
\begin{D}\label{carmes}
A positive measure $\mu$ on $\Di$ is a Carleson measure if
$\mu(s_{h})\leq Ch$ for every sector $s_{h}=\{re^{i\theta}\in\Di\ :\ 
1-h<r<1, |\theta-\theta_{0}|\leq h\}$. The least constant $C$ is called
the norm of Carleson measure $\mu$.
\end{D}
The following theorem of Carleson characterizes interpolating sequences,
see [C1].
\begin{Th}\label{carlchar}
The following are equivalent for a sequence of points $\{z_{j}\}$ in $\Di$:\\
(a)\ \ \ $\{z_{j}\}$ is an interpolating sequence;\\
(b)\ \ \ $\displaystyle \prod_{k,k\neq j}\left|\frac{z_{k}-z_{j}}
{1-\overline{z_{j}}z_{k}}\right|\geq\delta>0$\ for all $j$;\\
(c)\ \ \ $\displaystyle \left|\frac{z_{k}-z_{j}}
{1-\overline{z_{j}}z_{k}}\right|\geq a>0$ for all $k\neq j$ and the 
measure $\displaystyle \mu=\sum_{j=1}^{\infty}(1-|z_{j}|^{2})\delta_{z_{j}}$
is a Carleson measure. (Here $\delta_{z_{j}}$ is the Dirac measure at $z_{j}$.)
\end{Th}
The number $\delta$ in $(b)$ will be called the {\em characteristic of  
interpolating sequence} $\{z_{j}\}$.
\begin{Lm}\label{Gin}
For any $\xi\in G_{in}$ and $g\in G$ the element $\xi g$ belongs to $G_{in}$.
\end{Lm}
{\bf Proof.} Without loss of generality, assume that 
$\xi\in\beta G\setminus G$. Let $\{g_{j}\}_{j\in J}\subset G$ be a net
converging to $\xi$ such that $\{g_{j}(0)\}\subset\Di$ is an interpolating
sequence. Consider sequence $\{g_{j}g\}$, $g\neq e$. By definition,
$0<\rho(g_{j}(0),g_{j}g(0))=\rho(0,g(0)):=a<1$. Now according to
Corollary 1.6 of Ch. X and Lemma 5.3 of Ch. VII of [G], sequence
$\{g_{j}(0)\}$ can be decomposed into a finite disjoint union of subsets
$N_{i}$, i=1,...,k, such that every $N_{i}$ is an interpolating sequence
with the characteristic $\delta_{i}>\frac{2a}{1+a^{2}}$ and for any fixed $i$
the sequence $\{g_{l}g(0)\}$ with $g_{l}(0)\in N_{i}$ is interpolating. 
Clearly there is a number $s$ such that closure 
$\overline{\{g_{l}g(0)\}}\subset G_{in}$, $g_{l}\in N_{s}$, contains $\xi g$. 
Then, by definition, $\xi g\in G_{in}$.

The proof of the lemma is complete.\ \ \ \ \ $\Box$

For any interpolating sequence $\{z_{j}\}\subset\Di$ consider its closure 
$\overline{\{z_{j}\}}\subset E(M,\beta G)$. 
Let $I\subset E(M,\beta G)$ be union of all such closures.
\begin{Lm}\label{anal1}
Algebra $H^{\infty}\subset C(E(M,\beta G))$ separates points of $I$.
\end{Lm}
{\bf Proof.}
Let $\eta_{1},\eta_{2}\in I$ be different points. Since $E$ is a compact
Hausdorff space there are open neighbourhoods $V_{i}$ of $\eta_{i}$, $i=1,2$,
such that $\overline{V_{1}}\cap\overline{V_{2}}=\emptyset$. Without
loss of generality we may assume that $V_{i}\cong O_{i}\times H_{i}$,
where $O_{i}, H_{i}$ are open subsets of $M$ and $\beta G$, respectively
(see the description of topology on $E(M,\beta G)$).
Let now $x_{i}=\{x_{ji}\}_{j}\subset V_{i}\cap\Di$ be 
an interpolating sequence such that $\eta_{i}$ is its limit point, $i=1,2$.
We prove that $x_{1}\cup x_{2}$ is also an interpolating sequence. In fact, 
according
to property $(c)$ of Theorem \ref{carlchar} it suffices to check that
\begin{equation}\label{condit}
\left|\frac{x_{k1}-x_{j2}}
{1-\overline{x_{j2}}x_{k1}}\right|\geq a>0\ \ \ {\rm  for\ all}\ \ \
k,j.
\end{equation}
Let $\widetilde p:E(M,\beta G)\longrightarrow M$ be the projection. Consider 
two cases.\\
1. $\widetilde p(\eta_{1})\neq\widetilde p(\eta_{2})$.
In this case we can choose $V_{1}$ and $V_{2}$ such that
$\overline{O_{1}}\cap\overline{O_{2}}=\emptyset$, where 
$O_{i}:=\widetilde p(V_{i})$, i=1,2.
Assume, to the contrary,
that (\ref{condit})  is wrong. Then the pseudohyperbolic distance
$d(V_{1}\cap\Di,V_{2}\cap\Di):=
\inf_{z\in V_{1}\cap\Di,w\in V_{2}\cap\Di}\rho(z,w)$
between $V_{1}\cap\Di$ and $V_{2}\cap\Di$ equals 0. In particular,
$d(p^{-1}(O_{1}),p^{-1}(O_{2}))=0$, where $p=\widetilde p|_{\Di}$. 
Since $\rho$ is equivalent to the Euclidean distance on each compact
subset of $\Di$ and is invariant with respect to the action of the 
group of M\"{o}bius transformations, the above condition implies that
$\overline{p^{-1}(O_{1})}\cap\overline{p^{-1}(O_{2})}\neq\emptyset$.
It remains to note that 
$\overline{p^{-1}(O_{i})}\subseteq p^{-1}(\overline{O_{i}})$.
This gives a contradiction to the condition 
$\overline{O_{1}}\cap\overline{O_{2}}=\emptyset$.\\
2. $\widetilde p(\eta_{1})=\widetilde p(\eta_{2})$. In this case we can choose
$V_{1},V_{2}$ such that $O_{1}=O_{2}$,
but $H_{1}\cap G$ and $H_{2}\cap G$ are non intersecting subsets. Since
$G$ acts discretely, cocompactly on $\Di$ and $\rho$ is invariant with
respect to this action, we have $\rho(gx,x)\geq a>0$ for every $g\neq e$
and $x\in\Di$. Further, let us consider the fiber
$o=p^{-1}(\widetilde p(\eta_{1}))\subset\Di$. We may
assume without loss of generality
that $O_{1}$ is so small that the distance $d(\{y\},o)<a/4$ for each
$y\in V_{i}\cap\Di$, $i=1,2$. Let $x\in V_{1}\cap\Di$ and $y\in V_{2}\cap\Di$
and $x_{1},y_{1}\in o$ be nearest to $x,y$ 
points (the infimum for the distance is attained).
Applying the triangle inequality for $\rho$ we obtain
$$
\rho(x,y)\geq |\rho(x,x_{1})-\rho(x_{1},y)|\geq
\rho(x_{1},y_{1})-\rho(y_{1},y)-\rho(x,x_{1})>a/2.
$$
Here $x_{1}=gy_{1}$ for some $g\neq e$ (this follows from our assumption 
$(H_{1}\cap G)\cap (H_{2}\cap G)=\emptyset$).
So we proved that the sequences $x_{1}$ and $x_{2}$ are separated in
the pseudohyperbolic metric that is equivalent to (\ref{condit}).

Thus $x_{1}\cup x_{2}$ is an interpolating sequence.
Let now $f\in H^{\infty}$ be such that $f|_{x_{1}}=0$ and $f|_{x_{2}}=1$. Then 
clearly $f(\eta_{1})=0$ and $f(\eta_{2})=1$. This completes the proof of the
lemma.\ \ \ \ \ $\Box$

According to Hoffman's theorem, see, e.g. [G], Ch. X, Th. 2.5, every point of
$M_{a}$ is in the closure of an interpolating sequence. On the other
hand, the above lemma shows that quotient map
$F:E(M,\beta G)\longrightarrow M(H^{\infty})$ is one-to-one on $I$ and 
image $F(I)$ coincides with $M_{a}$. 
\begin{Lm}\label{belong}
$I=E(M,G_{in})$.
\end{Lm}
{\bf Proof.} Let $w\in I$ be a point from the
disk $i_{\eta}(\Di)$, $\eta\in X_{G}$. We will prove that $w\in E(M,G_{in})$.\\
Let $\{x_{j}\}_{j\in J}\subset\Di$ be a subnet of an interpolating sequence
converging to $w$. Consider $\{p(x_{j})\}\subset M$. Since 
$\widetilde p:E(M,\beta G)\longrightarrow M$ is continuous, 
$\lim_{j}p(x_{j})=\widetilde p(w)$. Let $\delta$ be the characteristic of 
$\{x_{j}\}_{j\in J}$ defined by Theorem \ref{carlchar} (b). Then we can
find an open neighbourhood $U\subset M$ of $\widetilde p(w)$ and a
subnet $S$ of $\{x_{j}\}_{j\in J}$ (we denote its elements by the same letters)
such that:\\
(1)\ $p^{-1}(U)\subset\Di$ is disjoint union of open sets homeomorphic to $U$;\\
(2)\ for each $x_{j}\in S$ there is $w_{j}\in p^{-1}(\widetilde p(w))$ such
that $\rho(x_{j},w_{j})<\delta/3$.\\
Then by Lemma 5.3 of Ch. VII of [G], $\{w_{j}\}$ is also an interpolating
sequence with the characteristic $\geq\delta/3$.
Moreover, condition $\lim_{j} p(x_{j})=\widetilde p(w)$ implies
$\lim_{j}\rho(x_{j},w_{j})=0$. In particular, for any bounded on $\Di$ function
$f$ uniformly continuous with respect to $\rho$,
$$
\lim_{j}|f(x_{j})-f(w_{j})|=0.
$$
But, by definition of $E(M,\beta G)$, $\lim_{j}f(x_{j})=f(w)$. Therefore 
$\lim_{j}f(w_{j})=f(w)$ and so $\lim_{j}w_{j}=w$. Further,
without loss of generality we may assume that $\{w_{j}\}$ consists of 
different points.
Consider the orbit $\{g(0)\}_{g\in G}$. By definition, there is a positive
number $r<1$ such that this orbit is an $r$-net in $\Di$ with respect to 
$\rho$. In particular, for every $w_{j}$ there is $g_{j}\in G$ so that 
$\rho(w_{j},g_{j}(0))\leq r$. Now according to Corollary 1.6 of
Ch. X and Lemma 5.3 of Ch. VII of
[G] the set $\{w_{j}\}$ can be represented as a finite disjoint union
of subsets $N_{i}$, $i=1,...,k$, such that every $N_{i}$ is an interpolating
sequence with the characteristic $\delta_{i}>\frac{2r}{1+r^{2}}$ and any
sequence $\{u_{l}\}$ satisfying
$\rho(u_{l},x_{l})\leq r$ for $x_{l}\in N_{i}$ is interpolating. Let
$s$ be such that closure $\overline{N_{s}}\subset E(M,\beta G)$ 
contains $w$. Then there is a subnet $W$ of the net $\{w_{j}\}_{j\in J}$
converging to $w$ which consists of elements of $N_{s}$.
Let $\{g_{l}\}\subset G$ be such that  
$\rho(g_{l}(0),w_{l})\leq r$ for $w_{l}\in W$. Then the sequence 
$\{g_{l}(0)\}$ is interpolating.  
Denote by $S$ the set of indexes of the $\{g_{l}\}$. Consider now
the map $D:S\times\Di\longrightarrow\Di$ given by
$$
D(l,z):=g_{l}(z).
$$
Clearly, for each bounded function $f$ on $\Di$ uniformly continuous 
with respect to $\rho$ the function $D^{*}f$ generates a uniformly 
bounded uniformly continuous with respect to $\rho$ family of
functions on $\Di$. Therefore according to Lemma \ref{main},
there is a continuous map
$\widetilde D:\beta S\times\Di\longrightarrow E(M,\beta G)$ which extends
$D$. Let $z_{j}$ be such that $D(z_{j})=w_{j}$ $(j=1,2,...)$. Then 
$z_{j}$ is in
the compact $\beta S\times B_{r}$, where $B_{r}$ is the closed ball with
respect to $\rho$ of radius $r$ centered at $0\in\Di$. Consider 
limit point
$z\in \beta S\times\Di$ of $\{z_{j}\}$ such that $\widetilde D(z)=w$.
Further, let $\xi\times\Di\subset\beta S\times\Di$ be the disk containing $z$.
Then, by definition, $\widetilde D(\xi\times\Di)$ coincides with some 
$i_{\psi}(\Di)$.
Therefore $w\in i_{\psi}(\Di)\cap i_{\eta}(\Di)$ and so 
$\psi=\eta$ (by Theorem \ref{proper}). Let $\{h_{l}\}\subset\{g_{l}\}$ be 
a subset such that its indexes generate a subnet converging to $\xi$.
Then $\{h_{l}(0)\}$ is interpolating and $\lim_{l}h_{l}(\Di)=i_{\eta}(\Di)$.
In particular, $i_{\eta}(\Di)\subset E(M,G_{in})$ and so $w\in E(M,G_{in})$.\\
We proved that $I\subset E(M,G_{in})$.
It remains to prove that $E(M,G_{in})\subset I$.\\
Let $y\in i_{\eta}(z)\subset E(M,G_{in})$. By definition, there is a net
$\{g_{j}\}_{j\in J}$ such that $\{g_{j}(0)\}$ is an interpolating sequence,
$\lim_{j}g_{j}=\xi\in G_{in}$ and $\xi$ represents
$\eta\in X_{G}$. In particular, $\lim_{j}g_{j}(\Di)=i_{\eta}(\Di)$
and there is $z\in\Di$ so that $\lim_{j}g_{j}(z)=y$. Clearly,
$\rho(g_{j}(0),g_{j}(z))=\rho(0,z):=a<1$. Then arguments similar to those
used in the proof of Lemma \ref{Gin} show that there is a subsequence 
$\{g_{l}\}_{l\in L}$ of $\{g_{j}\}$ such that $\{g_{l}(z)\}$ is an 
interpolating sequence whose closure contains $y$. Thus $y\in I$.

The lemma is proved.\ \ \ \ \ $\Box$

Let $H\subset\beta G$ be closure of an interpolating sequence. Then $H$ is
open by the definition of topology on $\beta G$. Further, $G_{in}$
is union of all such $H$ and therefore it is also open in $\beta G$. Now
the set $E(M,G_{in})$ is a bundle over $M$ with an open fibre, and, so it is an
open subset of $E(M,\beta G)$. In fact, $E(M,G_{in})$ can be covered by a
finite number of open sets homeomorphic to $U\times G_{in}$, where
$U\subset M$ is open.

To complete the proof of part (a) of the theorem we establish the following
\begin{Lm}\label{separ}
Let $z\in E(M,\beta G)$ be such that $F(z)\in M_{a}$. Then 
$z\in E(M,G_{in})$.
\end{Lm}
This together with the previous lemmas will show that $F|_{E(M,G_{in})}$
is a homeomorphism onto $M_{a}$.\\
{\bf Proof.}
Assume, to the contrary, that $z\not\in E(M,G_{in})$ and let
$\{z_{j}\}\subset\Di$ be a net converging to $z$.
Since $F(z)\in M_{a}$, there is a net $\{x_{j}\}_{j\in J}$ from elements
of an interpolating sequence which converges to $F(z)$ in $M(H^{\infty})$.
Let $B$ be an interpolating Blaschke product with zeros $S=\{x_{j}\}$,
and suppose that
$$
\inf_{j}(1-|x_{j}|^{2})|B'(x_{j})|\geq\delta>0.
$$
Note that $\delta$ coincides with the characteristic of $\{x_{j}\}$. 
Now according to Lemma 1.4 of Ch. X of [G],
there are $\lambda=\lambda(\delta)$, $0<\lambda<1$, and $r=r(\delta)$,
$0<r<1$, such that the set $B_{r}:=\{z\in\Di\ ; |B(z)|<r\}$ is the union
of pairwise disjoint domains $V_{n}$, $x_{n}\in V_{n}$, and
$$
V_{n}\subset\{z\in\Di\ ; \rho(z,x_{n})<\lambda\}.
$$
Consider open set $U=\{v\in M(H^{\infty})\ ; |B(v)|<r\}$. Since 
$B(F(z))=0$, $U$ is an open neighbourhood  of $F(z)$. Furthermore, by 
definition, the net
$\{z_{j}\}$ converges to $F(z)$ in $M(H^{\infty})$ (we assume that $F|_{\Di}$ 
is identity map). Therefore there is a subnet of 
$\{z_{j}\}$ containing in $U$ which converges to $F(z)$. (Without loss of 
generality we may assume that the net itself satisfies this property.)
Further, $U\cap\Di$ coincides with $\cup_{n}V_{n}$ and
so $\{z_{j}\}\subset\cup_{n}V_{n}$. 
Let $g_{n}$ be a M\"{o}bius
transformation such that $g_{n}(0)=x_{n}$ and $S$ denote the set of
indexes of $\{x_{n}\}$. Consider the map $h:S\times\Di\longrightarrow\Di$
defined by $h(l,z)=g_{l}(z)$. Then arguments of Lemma \ref{belong} show
that $h$ admits a continuous extension $\widetilde h:\beta S\times\Di
\longrightarrow E(M,\beta G)$. Further $\{z_{j}\}\subset\cup_{n}V_{n}$
and for each $z_{j}\in Z_{n}:=\{z_{j}\}\cap V_{n}$ there is 
$y_{j}\in n\times B_{r}$ such that 
$h(y_{j})=z_{j}$. Here $B_{r}$ is the closed ball in $\Di$ with respect to
$\rho$ of radius $r$ centered at $0$.
Let $Y\subset\beta S\times B_{r}$ be closure of 
$\{y_{j}\}$ in $\beta S\times\Di$. Then $\widetilde h(Y)$ is a compact
in $E(M,\beta G)$ containing  $\{z_{j}\}$. In particular, there is
$y\in Y$ such that $\widetilde h(y)=z$. Assume that $y\in\xi\times\Di$,
$\xi\in\beta S$. Then the set $\widetilde h(\xi\times\Di)$ belongs to some 
disk $i_{\eta}(\Di)$, $\eta\in X_{G}$. But by definition, 
$\widetilde h(\xi\times 0)$ is a limit point of the interpolating sequence
$\{x_{j}\}$. Therefore $h(\xi\times 0)\in E(M,G_{in})$ and
$i_{\eta}(\Di)\subset E(M,G_{in})$. This shows that $z\in E(M,G_{in})$ and
contradicts to the original assumption.

The lemma is proved.\ \ \ \ \ $\Box$

So we proved that $F:E(M,G_{in})\longrightarrow M_{a}$ is one-to-one
and $F^{-1}(M_{a})=E(M,G_{in})$ is open in $E(M,\beta G)$. Therefore
$M_{a}\subset M(H^{\infty})$ is open and $F:E(M,G_{in})\longrightarrow
M_{a}$ is a homeomorphism. The proof of part (a) is complete.\\
Since $M_{a}\subset E(M,\beta G)$, statements (b) and (c) follow 
directly from the corresponding  statements of Theorem \ref{te1} and
Proposition \ref{proper} (c).

The proof of Theorem \ref{class} is complete.\ \ \ \ \ $\Box$
\begin{R}\label{theodim}
{\rm Based on Theorem \ref{class} one obtains an alternative proof of the 
result
of Su\'{a}rez asserting that $dim\ M(H^{\infty})=2$.  In fact,
according to Proposition \ref{proper} (a), $dim\ K\leq 2$ for any compact
$K\subset M_{a}$. Moreover, by Theorem 3.4 of [S1],
$dim\ M(H^{\infty})\setminus M_{a}=0$. Then 
$dim\ M(H^{\infty})=2$ by [N], Ch. 2, Th. 9-11.}
\end{R}
We complete this section by the result similar to Proposition \ref{proper}
(d).
\begin{Proposition}\label{hpath}
For any different points $x,y\in M(H^{\infty})$ there is a connected
one-dimensional compact set $l_{x,y}\subset M(H^{\infty})$ containing
them.
\end{Proposition}
{\bf Proof.}
Let $x_{1}, y_{1}\in E(M,\beta G)$ be such that $F(x_{1})=x$ and
$F(y_{1})=y$. By Proposition \ref{proper} (d), there is a connected
one-dimensional compact set $l\subset E(M,\beta G)$ containing
$x_{1},x_{2}$. Then $l_{x,y}:=F(l)\subset M(H^{\infty})$ satisfies the
required property. Actually, by Theorem \ref{class}, $dim\ K\leq 1$ for any 
compact $K\subset M_{a}\cap l_{x,y}$. But 
$dim\ l_{x,y}\cap (M(H^{\infty})\setminus M_{a})\leq 0$ by [S1], Th. 3.4.
Therefore $dim\ l_{x,y}\leq 1$ by [N], Ch. 2, Th. 9-11. Moreover,
$l_{x,y}$ is connected because $F$ is continuous. Thus $dim\ l_{x,y}=1$
(otherwise it is totally disconnected).\ \ \ \ \ $\Box$
%===============================
\sect{\hspace*{-1em}. Subsets of $M_{a}$.}
In this section we prove results of Section 2.2.\\
{\bf Proof of Theorem \ref{subsets}.} Let $X$ be a complex Riemann surface 
which admits a non-constant bounded holomorphic function. Let 
$r:\Di\longrightarrow X$ be the universal 
covering over $X$ and $H=\pi_{1}(X)$. We begin with the following
\begin{Lm}\label{surf}
For every $z\in X$ the sequence $\{r^{-1}(z)\}\subset\Di$ is interpolating.
\end{Lm}
{\bf Proof.}
Let $f$ be a non-constant bounded holomorphic function on $X$.
Then clearly $r^{*}f\in H^{\infty}$. Assume, without loss of
generality, that $r^{*}f$ vanishes at points of $r^{-1}(z)$. Further, we can 
also assume that $0\in r^{-1}(z)$. In fact, let $y\in r^{-1}(z)$ and $y\neq 0$.
Then there is a M\"{o}bius transformation $g$ of $\Di$, such that 
$g(y)=0$. Moreover, $g$ is equivariant with respect to the actions on $\Di$ of
groups $H$ and $gHg^{-1}$, respectively. In particular, it determines a
biholomorphic mapping $\widetilde g$ between quotient spaces $X=\Di/H$ and
$Y=\Di/gHg^{-1}$. Thus $Y$ satisfies conditions of our lemma. Now if we
prove that $\{ghg^{-1}(0)\}_{h\in H}$ is interpolating, then  
$r^{-1}(z)=\{h(y)\}_{h\in H}$ will be also interpolating. 

So let $r^{*}f$ be vanishing at $\{h(0)\}_{h\in H}$. Then according to Th.2.1, 
Ch. II of [G],
$\sum_{h\in H}(1-|h(0)|)<\infty .$
To prove that $\{h(0)\}_{h\in H}$ is
interpolating we check property (b) of Theorem \ref{carlchar}. Since
the pseudohyperbolic distance $\rho$ is invariant with respect to the action
of the group of M\"{o}bius transformations, it suffices to check that
$\prod_{h\neq e}|h(0)|\neq 0$. Let $H_{1}\subset H$ be such that
$1-|h(0)|<1/2$ for any $h\in H_{1}$. Set $H_{2}=H\setminus H_{1}$.
Note that if $h\in H_{1}$ then $\log|h(0)|>-2(1-|h(0)|)$. Summing
these inequalities we obtain
$$
\sum_{h\in H_{1}}\log |h(0)|>-2\sum_{h\in H_{1}}(1-|h(0)|)>-\infty.
$$
Taking the exponent in the above inequality and noting that $H_{2}$ is a 
finite set we get
$$
\prod_{h\neq e}|h(0)|\neq 0.
$$
Therefore $\{h(0)\}$ is interpolating. According to the remark given in the
proof this case is general so $\{h(y)\}$ is interpolating for any $y\in\Di$. 

The lemma is proved.\ \ \ \ \ $\Box$
\begin{R}\label{after}
{\rm Similarly to the proof of case (1) of Lemma \ref{anal1} one 
deduces from the above lemma that for any $z_{1},z_{2}\in X$ sequence
$r^{-1}(z_{1})\cup r^{-1}(z_{2})$ is interpolating.}
\end{R}
Consider now the bundle $E(X,\beta H)$ over $X$ with fibre $\beta H$.
Then as before we have\\
(a)\ $\Di$ is embedded into $E(X,\beta H)$ as an open everywhere dense subset;
\\
(b)\ every $f\in H^{\infty}$ has a continuous extension to $E(X,\beta H)$
(denoted by the same letter)
such that for each $\eta\in X_{H}$, the function $f\circ i_{\eta}$ belongs to
$H^{\infty}$;\\
(c)\ $dim\ E(X,\beta H)=2$.\\
According to Lemma \ref{surf} and Remark \ref{after}
$H^{\infty}$ separates points of 
$E(X,\beta H)$. Therefore there is a continuous injective mapping 
$I:E(X,\beta H)\longrightarrow M_{a}$. By Theorem \ref{class} (c)
for any $\eta\in\beta H/H$, there is $m\in M(H^{\infty})$ such that
$I(i_{\eta}(\Di))$ belongs to Gleason part $P(m)$. Moreover, let
$\xi\in\beta H$ represent $\eta$ and $\{h_{j}\}_{j\in J}$ be a net of
elements of $H$ converging to $\xi$. By definition, 
$I(i_{\eta}(z))=\lim_{j}h_{j}(z)$ in $M(H^{\infty})$ for any $z\in\Di$. 
Then the arguments of Th.1.7 of Ch. X of [G] show that $I(i_{\eta}(\Di))=P(m)$.
Obviously, for any $f\in H^{\infty}$ the function 
$f\circ I\circ i_{\eta}\in H^{\infty}$.
It remains to prove that 
for any open $U\subset E(X,\beta H)$ its image $I(U)$ is open 
in $M_{a}$.

It suffices to prove the statement for open sets generating a basis
of topology on $E(X,\beta H)$. We now describe such sets.
Let $\widetilde r:E(X,\beta H)\longrightarrow X$ be the projection and
$r=\widetilde r|_{\Di}$. Consider a point $z\in X$ and the set
$r^{-1}(z)=\{z_{i}\}\subset\Di$. Let $\epsilon$ be so small
that for $V_{i}:=\{z\in\Di\ ; \rho(z,z_{i})<\epsilon\}$ we have
$V_{i}\cap V_{j}=\emptyset$ for $i\neq j$ and for each $i$ the mapping 
$r|_{V_{i}}$ maps univalently $V_{i}$ 
onto an open neighbourhood $U$ of $z\in X$. 
Let $\{z_{n_{i}}\}_{i\geq 1}$ be a subset of $\{z_{n}\}$.
Denote $V=\cup_{i}V_{n_{i}}$ and consider closure $\overline V$ in 
$E(X,\beta H)$. Finally set 
$V'=\widetilde r^{-1}(U)\cap\overline{V}$. Then $V'\subset E(X,\beta H)$ is 
open and the sets of this form determine a basis of topology on 
$E(X,\beta H)$. So it suffices to prove that $I(V')$ is open for any such 
$V'$. Let $y\in I(V')$
and $y'=\widetilde r(I^{-1}(y))\in U\subset X$. Denote by $B_{y'}$
the interpolating Blaschke product with zeros at 
$\{y_{i}\}=\{r^{-1}(y')\cap V_{n_{i}}\}\subset\Di$.
Then clearly $y$ is a limit point of $\{y_{i}\}$ and the set
$Y=\{v\in M(H^{\infty}); |B_{y'}(v)|<\epsilon'\}$ is an open neighbourhood of
$y$. We will prove that if $\epsilon'>0$ is sufficiently small then 
$Y\subset I(V')$. In fact,
according to Lemma 1.4 of Ch. X of [G], for sufficiently small 
$\epsilon'>0$ intersection $Y\cap\Di$ is 
disjoint union of sets $Y_{n}\subset\{z\in\Di\ ; \rho(y_{n},z)<\lambda\}$,
where $\lambda=\lambda (\epsilon')>0$.
Moreover, if $\epsilon'\to 0$ then we can also choose  $\lambda\to 0$.
So we can take $\epsilon'$ so small that $Y\cap\Di\subset V_{1}\subset V'$,
where $V_{1}$ is a compact subset of $V'$. 
Therefore $\overline{Y\cap\Di}$ contains $Y$ and belongs to the compact
$I(V_{1})$. In particular, $Y\subset I(V')$. This shows that $I(V')$ is open. 

Theorem \ref{subsets} is proved.\ \ \ \ \ $\Box$
\begin{R}
{\rm Using the arguments of Section 5 we can point out embedded and 
non-embedded analytic disks in $E(X,\beta H)\subset M(H^{\infty})$.}
\end{R}
{\bf Proof of Corollary \ref{free}.}
Let $X\subset\Di$ be an open set obtained by removing from $\Di$ a finite
or countable number of disjoint closed disks. Then clearly $X$ admits a
non-constant holomorphic function and $H=\pi_{1}(X)$ is free 
group with finite or countable number of generators. By Theorem
\ref{subsets} $E(X,\beta H)$ is an open subset of $M(H^{\infty})$. Let
now $\xi_{1}$ and $\xi_{2}$ be different points of $\beta H$ such that
$\overline{P(\xi_{1})}\cap\overline{P(\xi_{2})}\neq\emptyset$, where
$P(\xi_{i}):=\{\xi_{i}h,\ h\in H\}$, $i=1,2$. We have to prove that
one of these sets contains the other. Let $\widetilde\xi_{i}$ be image
of $\xi$ in the quotient set $X_{H}=\beta H/H$, $i=1,2$.
Consider analytic disks
$P_{i}:=i_{\widetilde\xi_{i}}(\Di)\subset E(X,\beta H)\subset M_{a}$, $i=1,2$.
Then every such disk is a non-trivial Gleason part. According to our 
assumption closures $\overline{P_{1}},\overline{P_{2}}$ in $M(H^{\infty})$ 
have non-empty intersection. 
Therefore the above assumption together with Corollary 2.7 of [S2] imply
that one of these sets contains the other.
But since $E(X,\beta H)$ is an open subset of $M(H^{\infty})$, we
have $\overline{P_{i}}\cap E(X,\beta H)$, $i=1,2$, coincides with
closure $P_{i}'$ of $P_{i}$ in $E(X,\beta H)$. In particular,
one of the sets $P_{1}',P_{2}'$ contains the other. Assume, e.g., that 
$P_{2}'\subset P_{1}'$. Further, let $y\in X$ and $U_{y}\subset X$ be an open 
disk centered at $y$. Then bundle $E(X,\beta H)$ over $U_{y}$ is homeomorphic
to $U_{y}\times\beta H$ (without loss of generality we may identify these
sets). According to the above arguments, $\{(y,\xi_{2}h)\}_{h\in H}$
is in the 
closure of $\{U_{y}\times\xi_{1}g\}_{g\in H}$. Therefore for any
$h\in H$ there is a net $\{(z_{j},\xi_{1}g_{j})\}_{j\in J}$, $z_{j}\in U_{y},
g_{j}\in H$, converging to $(y,\xi_{2}h)$. From here it follows immediately
that $\lim_{j}z_{j}=y$. Thus $\lim_{j}(y,\xi_{1}g_{j})=(y,\xi_{2}h)$.
This shows that $\overline{P(\xi_{2})}\subset\overline{P(\xi_{1})}$.

The corollary is proved.\ \ \ \ \ $\Box$\\
{\bf Proof of Theorem \ref{algeb}.}
First prove that $\Di$ is dense in the maximal ideal space $M(A_{Y})$ of
the algebra $A_{Y}$. It is equivalent to the following problem:\\
Let $f_{1},...f_{n}$ be holomorphic functions from $A_{Y}$ satisfying
\begin{equation}\label{alg1}
\max_{i}|f_{i}(z)|>\delta\ \ \ \ (z\in\Di).
\end{equation}
Find functions $g_{1},...,g_{n}$ in $A_{Y}$ so that
$$
\sum_{i}f_{i}g_{i}=1.
$$
Since $A_{Y}$ is closure in $C(\Di)$ of the algebra generated by 
all $a^{*}(H_{U})$, it suffices to solve the above problem under assumption:
there is a domain $U$, 
$\overline{Y}\subset U\subset Z$, 
and bounded holomorphic functions $f_{1}',...,f_{n}'$ defined on 
$r_{Z}^{-1}(U)$ and satisfying (\ref{alg1}) there such that 
$a^{*}(f_{i}')=f_{i}$, $i=1,...,n$. In particular, it suffices to find a
domain $V$, $\overline{Y}\subset V\subset U$, and bounded holomorphic 
functions $g_{1}',...,g_{n}'$ defined on $r_{Z}^{-1}(V)$ such that
$\sum g_{i}'f_{i}'=1$.
Further, similarly to constructions of Section 6, define a Banach holomorphic
vector bundle $E_{U}$ over $U$ with fiber $l^{\infty}(H)$ associated with 
the universal covering $r_{Z}:r_{Z}^{-1}(U)\longrightarrow U$. Then
we can think of a bounded holomorphic function on $r_{Z}^{-1}(U)$ as a 
holomorphic section of $E_{U}$. Moreover, the  natural
multiplication defined on $l^{\infty}$ determines a multiplication on
sections of $E_{U}$. In what follows we identify bounded functions on
$r_{Z}^{-1}(U)$ with sections of $E_{U}$.
Set
$\phi_{i}(z)=\overline{f_{i}'}(z)/\sum\ |f_{k}'(z)|^{2},\ 
i=1,...,n$. Then according to (\ref{alg1}) each $\phi_{i}$ is a smooth
section of $E_{U}$ (existence of derivatives follows from the
normality) and $\sum\phi_{i}f_{i}'=1$. In particular, 
$\overline{\partial}\phi_{i}$ is a well-defined $(0,1)$-form with values in 
$E_{U}$.  Consider now $E_{U}$-valued $\overline{\partial}$-equations
\begin{equation}\label{coro}
\overline{\partial}b_{j,k}=\phi_{j}\overline{\partial}\phi_{k},\ \ \ \ 
1\leq j,k\leq n.
\end{equation}
Since $U$ is a Stein manifold, the results of
[B], sect. 4 imply that there are smooth sections $b_{j,k}$ of $E_{U}$ 
(not necessary bounded) that solve (\ref{coro}).
(This is equivalent to the identity
$H^{1}(U,{\cal O}_{E_{U}})=0$, where ${\cal O}_{E_{U}}$ is the sheaf of germs
of local holomorphic sections of $E_{U}$.)
Let now a domain $V$ be such that
$\overline{Y}\subset\overline{V}\subset U$. Then $b_{j,k}|_{V}$
is bounded ($b_{j,k}(V)\subset E_{U}$ is a precompact). Therefore 
$$
g_{j}'(z)=\phi_{j}(z)+\sum_{k=1}^{n}(b_{j,k}-b_{k,j})f_{k}(z)',\ \ \
j=1,...,n,
$$
are bounded holomorphic sections over $V$ satisfying 
$\sum\ (g_{j}'f_{j}')|_{V}=1$ (see also [G], Ch. VIII).
Further, according to our identification $g_{j}'$ represents a bounded
holomorphic function on $r_{Z}^{-1}(V)$.  Thus the functions 
$a^{*}(g_{j}')\in A_{Y}$, $j=1,...,n$, solve the required corona problem.

So we proved that $\Di$ is dense in $M(A_{Y})$. Note also that 
any function $f\in A_{Y}$ admits a continuous extension to 
$E(\overline{Y},\beta H)$ and according to Lemma \ref{surf} and Remark 
\ref{after} the extended algebra separates points of 
$E(\overline{Y},\beta H)$.
Since the above space is a compact and $\Di$ is dense in $M(A_{Y})$,  
the space $E(\overline{Y},\beta H)$ is homeomorphic to the maximal 
ideal space $M(A_{Y})$. Then Proposition \ref{proper} (a) asserts
that $dim\ M(A_{Y})=2$. Therefore Theorem 3 of [L] and the results of
Section 3 show that the problem of the complementation of matrices is valid
for $A_{Y}$.

The proof is complete.\ \ \ \ \ $\Box$
\begin{R}\label{merg}
{\rm 
Let $\widetilde r_{Z}:E(\overline{Y},\beta H)\longrightarrow\overline{Y}$ be 
the map extending $r_{Z}$. Then the \v{S}ilov boundary of $A_{Y}$ coincides 
with $\widetilde r_{Z}^{-1}(\overline{Y}\setminus Y)$.\\
Assume that $Y$ is a domain in $\Co$ satisfying conditions of
Theorem \ref{algeb} and $\Co\setminus\overline{Y}$ is a disjoint union of
connected domains. Then using some approximation theorems one can
show that the extension of $A_{Y}$ to $E(\overline{Y},\beta H)$ coincides
with the algebra of holomorphic functions continuous up to the boundary
$\widetilde r_{Z}^{-1}(\overline{Y}\setminus Y)$.}
\end{R}
%==============================================
%==========================

\noindent Department of Mathematics\\
Ben Gurion University of the Negev, P.O.B. 653,\\
Beer-Sheva 84105, Israel.
\end{document}